\documentclass[reqno]{amsart}

\usepackage[utf8]{inputenc}
\usepackage[T1]{fontenc}
\usepackage[english]{babel}
\usepackage[backend=biber, style=numeric, sorting=nty, maxbibnames=10, doi=false,isbn=false,url=false,eprint=false, giveninits=true]{biblatex}
\usepackage{csquotes}
\usepackage{latexsym}
\usepackage{amssymb}
\usepackage{dsfont}
\usepackage{mathtools}
\usepackage[hidelinks]{hyperref}
\usepackage{enumitem}
\usepackage{linegoal}
\usepackage{float}
\usepackage{multicol}
\usepackage{mathrsfs} 
\usepackage{stmaryrd} 
\usepackage{tikz}
\usetikzlibrary{patterns}
\usetikzlibrary{calc}
\usetikzlibrary{decorations.markings}
\usepackage{tikz-3dplot}
\usetikzlibrary{3d, arrows.meta, calc}

\hypersetup{final} 

\DeclareFieldFormat[unpublished]{title}{\emph{#1}}
\AtEveryBibitem{
	\ifentrytype{article}{}{\clearfield{pages}}
	\clearfield{edition}
	\clearfield{series}
	\ifentrytype{unpublished}{}{\clearfield{note}}
	\ifentrytype{book}{\clearfield{volume}}{}
	\ifentrytype{unpublished}{}{\clearfield{url}}
}
\allowdisplaybreaks[4]

\theoremstyle{plain}
\newtheorem{theorem}{Theorem}[section]
\newtheorem*{proposition*}{Proposition}
\newtheorem{proposition}[theorem]{Proposition}
\newtheorem{lemma}[theorem]{Lemma}

\newtheorem{definition}[theorem]{Definition} 
\theoremstyle{definition}
\newtheorem{remark}[theorem]{Remark}	   

\numberwithin{equation}{section} 

\newcommand{\R}{\mathbb{R}}

\newcommand{\N}{\mathbb{N}}

\newcommand{\Rsym}{\R^{3\times3}_{\mathrm{sym}}}

\newcommand{\cC}{\mathcal{C}}
\newcommand{\cD}{\mathcal{D}}

\newcommand{\cG}{\mathcal{G}}

\newcommand{\cL}{\mathcal{L}}
\newcommand{\cM}{\mathcal{M}}
\newcommand{\cN}{\mathcal{N}}

\renewcommand{\div}{\operatorname{div}}
\DeclareMathOperator{\curl}{curl}
\newcommand{\tr}{\operatorname{tr}}
\newcommand{\dist}{\operatorname{dist}}

\newcommand*\dd{\mathop{}\!\mathrm{d}}

\newcommand{\pt}{\partial_t}

\newcommand{\p}{\partial}

\newcommand{\e}{\mathrm{e}}

\let\latexepsilon\epsilon
\let\latexvarepsilon\varepsilon
\renewcommand{\epsilon}{\latexvarepsilon}
\renewcommand{\varepsilon}{\latexepsilon}

\renewcommand{\complement}{\mathsf{c}} 
\newcommand{\supp}{\operatorname{supp}}

\let\oldexists\exists \let\exists\relax \DeclareMathOperator{\exists}{\oldexists} 
\let\oldforall\forall \let\forall\relax \DeclareMathOperator{\forall}{\oldforall}

\newcommand{\inv}[1]{{#1}^{-1}}

\newcommand{\ol}{\overline}

\let\latexchi\chi
\makeatletter
\renewcommand\chi{\@ifnextchar_\sub@chi{{\latexchi}}}
\newcommand{\sub@chi}[2]{
	{\@ifnextchar^{\subsup@chi{#2}}{{\latexchi}^{}_{#2}}}
}
\newcommand{\subsup@chi}[3]{
	{{\latexchi}_{#1}^{#3}}
}

\makeatletter
\def\moverlay{\mathpalette\mov@rlay}
\def\mov@rlay#1#2{\leavevmode\vtop{%
		\baselineskip\z@skip \lineskiplimit-\maxdimen
		\ialign{\hfil$\m@th#1##$\hfil\cr#2\crcr}}}
\newcommand{\charfusion}[3][\mathord]{
	#1{\ifx#1\mathop\vphantom{#2}\fi
		\mathpalette\mov@rlay{#2\cr#3}
	}
	\ifx#1\mathop\expandafter\displaylimits\fi}
\makeatother
\newcommand{\cupdot}{\charfusion[\mathbin]{\cup}{\cdot}}

\newcommand{\abs}[1]{ \left\lvert#1\right\rvert} 
\newcommand{\norm}[1]{\left\lVert#1\right\rVert} 
\newcommand{\smallnorm}[1]{\lVert#1\rVert} 
\newcommand{\bignorm}[1]{\big\lVert#1\big\rVert} 
\newcommand{\biggnorm}[1]{\bigg\lVert#1\bigg\rVert} 
\newcommand{\Bignorm}[1]{\Big\lVert#1\Big\rVert} 
\newcommand{\vertiii}[1]{{\left\vert\kern-0.25ex\left\vert\kern-0.25ex\left\vert #1 
		\right\vert\kern-0.25ex\right\vert\kern-0.25ex\right\vert}}

\AtBeginDocument{
}

\newcommand\restr[2]{{
		\left.\kern-\nulldelimiterspace 
		#1 
		\vphantom{\big|} 
		\right|_{#2} 
}}

\newcommand{\seq}{\hphantom{={}}}

\newcommand{\alphad}{\alpha_{\mathrm{d}}}
\newcommand{\betad}{\beta_{\mathrm{d}}}

\newcommand{\bndInh}{h}

\begin{document}

\title[Approximation of the Nonautonomous Maxwell Interface Problem]{Approximation of the Nonautonomous Maxwell System with an Interface}
\author{Richard Nutt}
\address{Richard Nutt\\
	Karlsruhe Institute of Technology\\
	Department of Mathematics\\
	Englerstraße~2\\
	76131 Karlsruhe\\
	Germany}
\email{richard.nutt@kit.edu}




\begin{abstract}
	To obtain higher-regularity solutions for nonautonomous Maxwell systems with an interface, we approximate the initial data and coefficients by smoother counterparts. To preserve compatibility conditions of the approximations, required to obtain solutions to the Maxwell systems, we systematically adjust the approximating sequence of initial data. While the existing well-posedness theory in \cite{schnaubelt2022} guarantees more regular solutions, standard a-priori estimates fail to ensure convergence in the desired norm. We overcome this by a mollification argument and a div-curl estimate from \cite{bresch2026} and \cite{dohnal2024}.
\end{abstract}

\maketitle

\section{Introduction}

The Maxwell system is the foundation of electromagnetic theory.
It poses significant challenges for the analytical treatment, e.g.\ it is a highly coupled system with difficult boundary or interface conditions.
To obtain solutions in higher regularity the data thus has to satisfy sophisticated compatibility conditions, see Definition~\ref{def:CC}.
As in other hyperbolic problems, the derivation of the energy equality or other properties of solutions requires more regularity than one expects or needs for the respective formulas.
Thus, one has to regularize the data, and possibly the coefficients, and then show that the resulting solutions of higher regularity converge to the given one in the right norm.
Here, the compatibility conditions are a serious obstacle.

For systems with linear, autonomous coefficients one can often exploit the density of the domains $D(A^n)$ of the generator $A$ of the system.
However, this approach breaks down for nonautonomous problems like the interface system
\begin{alignat}{-1}
	\begin{split} \label{eq:basic-maxwell}
		\pt \left(\epsilon(t,x) E(t,x)\right) & = \curl H(t,x) - \sigma(t,x)E(t,x) - f(t,x), \\
		\pt\left(\mu(t,x) H(t,x)\right)       & = - \curl E(t,x)- g(t,x),
	\end{split}
	\hspace{0.5cm}                                                                                                                              &                                           & t\geq 0, x \in \Omega_\pm,                                   \\
	\notag E(t,x) \times \nu(x)                                                                                                                 & = \bndInh(t,x),                           &                            & t\geq 0, x \in \partial \Omega, \\
	\begin{split}\notag
		\llbracket E(t,x) \times \nu^\Gamma(x) \rrbracket & =\kappa(t,x), \quad\llbracket H(t,x) \times \nu^\Gamma(x) \rrbracket = \varkappa(t,x),
	\end{split} &                                           & t\geq 0, x \in \Gamma,                                                                                                             \\
	\notag E(0,x)                                                                                                                               & = E_0(x), \hspace{19.5mm} H(0,x)= H_0(x), &                            & x \in \Omega_\pm,
\end{alignat}
with instantaneous material laws.
Such systems in particular occur as linearizations of quasilinear problems with permittivity $\epsilon(t,x) = \epsilon(x,E(t,x))$ and permeability $\mu(t,x) = \mu(x,H(t,x))$, see \cite{lasiecka2019} and \cite{nutt2026-2}, where the coefficients depend on the fields.
In this setting approximation results are so far restricted to very simple cases, see e.g.\ the proof of Lemma~4.1 in \cite{lasiecka2019} or Lemma~5.8 in \cite{schnaubelt2022}.

In this paper we provide a rather complete discussion of the problem in the context of \eqref{eq:basic-maxwell}, see also Remark~\ref{rem:quasilinear} for the application to quasilinear problems.

In order to obtain solutions of higher regularity we will approximate the data and coefficients by more regular ones.
A major obstacle is that the approximation must also satisfy the compatibility conditions of Definition~\ref{def:CC} depending on the higher regularity level.
We therefore iteratively adjust the initial data of the approximating sequence to satisfy the compatibility conditions in Proposition~\ref{prop:fix-CC} by adding small correction terms to the approximation.
Here we heavily use the structure of the curl operator.
Up to the order $m$ of the given solution, the correction terms can be chosen small by exploiting that the approximation only perturbs the compatibility conditions by an arbitrarily small amount.
For orders greater than $m$, this approach no longer works.
Instead, we rely on a scaling argument to control the higher-order corrections.

Then, the well-posedness theory from Theorem~3.1 in \cite{schnaubelt2022} yields solutions with additional regularity, but the a-priori estimate given there does not guarantee convergence in the right norm, as outlined at the beginning of Section~\ref{sec:3-convergence-result}.
Therefore, to establish the convergence, we first use a mollification and Gronwall argument to show the convergence of the highest-order time derivatives.
We then extend this result to the full space-time-norm by employing a div-curl estimate from \cite{bresch2026} and \cite{dohnal2024}. In order to apply the estimate we have to first add a correction term due to the inhomogeneities on the boundary and interface.
In an intricate induction argument we then trade spacial against temporal regularity by employing the structure of \eqref{eq:basic-maxwell}.
Moreover, terms of highest order are once more controlled by a Gronwall argument, finally allowing us to show convergence in the desired norm.

At the end we discuss the application of the approximation result to the quasilinear setting, where one usually assumes regularity on $\partial_\xi \epsilon(x,\xi) $ and not $\partial_t \epsilon(x,E(t,x))$ (respectively for $\mu$) as we do in our discussion of the nonautonomous, linear system.

Our approximation result is used in the companion contribution \cite{nutt2026-2} where we show that solutions of the quasilinear system decay exponentially to $0$ for small data in the presence of a conductivity near the boundary.

\section{Problem Setting and Compatibility Conditions}
Let $m, k \in \N_0$ with $m \leq k$. We study a bounded open set $\Omega \subseteq \R^3$ which is split by an interface $\Gamma$ into disjoint, possibly disconnected sets $\upsilon$ and $\omega$ with $\partial \Omega \subseteq \partial \omega$. We set $\Omega_\pm \coloneqq \upsilon \cupdot \omega = \Omega \setminus \Gamma$.
It is assumed that $\partial \Omega \cap \Gamma = \emptyset$ and $\omega$ to have a $C^{k+1}$-boundary.
Let $\nu$ and $\nu^\Gamma$ denote the outer unit normal of $\Omega$ and $\upsilon$ respectively.

\vspace*{-7mm} 
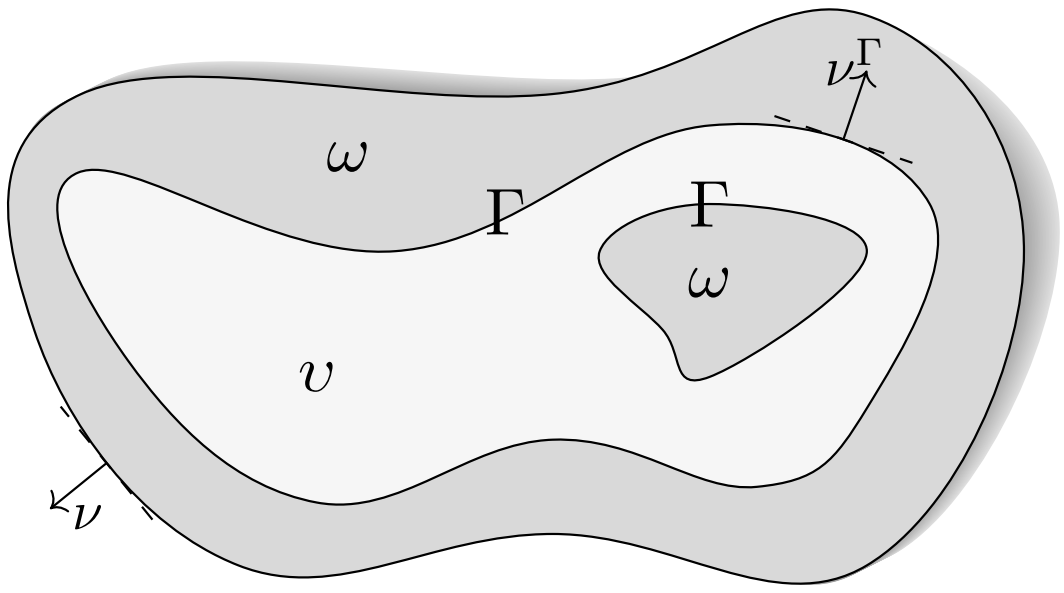
\begin{figure}[H]
	\centering
	\resizebox{0.6\textwidth}{!}{%
		\def\nodeScale{1.2}
		\newcommand\decWidth{8.45pt}

		\begin{tikzpicture}[decoration={markings,
						mark= between positions 0 and 1 step 2pt with{
								\path[pattern=north west lines, pattern color=black!30] (0,-\decWidth) rectangle (3pt,\decWidth);
							}}]


			\foreach \step in {25, ..., 0} {
					\pgfmathsetmacro{\xoff}{\step*0.010+0.28}
					\pgfmathsetmacro{\yoff}{\step*0.005+0.14}
					\pgfmathtruncatemacro{\col}{75 - \step*2}
					\fill[scale=0.88, gray!\col] plot [smooth cycle,tension=0.7, scale=1.12] coordinates {
							(\xoff+2,\yoff+1) (\xoff+4,\yoff+1.2) (\xoff+5.78,\yoff+0.85) (\xoff+6.8,\yoff+3)
							(\xoff+5.5,\yoff+4.31) (\xoff+4,\yoff+3.9) (\xoff+0.69,\yoff+3.9) (\xoff+0.7,\yoff+2.5)
						};
				}


			\draw[fill=gray!30] plot [smooth cycle,tension=0.7] coordinates {(2,1) (4,1.2) (6,1) (7,3) (6,4.5) (4,4) (1,4) (0.7,2.5)};

			\draw[fill=gray!7] plot [smooth cycle,tension=0.7] coordinates {(2.5,1.4) (4,1.8) (5.3,1.5) (6,2) (6.4,3.3) (5,3.8) (3,3) (0.98,3.5) (1.3,2.3)};

			\draw[fill=gray!30] plot [smooth cycle,tension=0.7] coordinates {(4.7,2.5) (5,2.2) (6,3) (5,3.3) (4.3,3) };


			\coordinate (A) at (5.851,3.711);
			\coordinate (B) at (6,4.15);
			\draw[->] (A) -- node[above] {$\nu^\Gamma$} (B);
			\draw[dashed] ($(A)!1!90:(B)$) -- ($(A)!1!-90:(B)$);

			\coordinate (C) at (1.161,1.651);
			\path let \p1 = ($(B)-(A)$), \n1 = {veclen(\x1,\y1)}
			in coordinate (D) at ($(C)!\n1!(0.8,1.358)$);

			\draw[->] (C) -- node[below] {$\ \nu$}(D);
			\draw[dashed] ($(C)!1!90:(D)$) -- ($(C)!1!-90:(D)$);

			\node[scale=\nodeScale] at (2.5,2.2) {$\upsilon$};
			\node[scale=\nodeScale] at (3.7,3.25) {$\Gamma$};
			\node[scale=\nodeScale] at (5,3.3) {$\Gamma$};
			\node[scale=\nodeScale] at (2.7,3.6) {$\omega$};
			\node[scale=\nodeScale] at (5,2.8) {$\omega$};

		\end{tikzpicture}
	}
	\vspace*{-2mm}
	\caption{$2d$-cut through $\Omega \subset \mathbb{R}^3$.}
\end{figure}
In view of \eqref{eq:basic-maxwell} we consider the nonautonomous, anisotropic Maxwell system
\begin{alignat}{-1}
	\begin{split}\label{eq:maxwell-original}
		\pt \left(\alpha(t,x) u(t,x)\right) & = \hphantom{-}\curl v(t,x) - \tilde \sigma(t,x)u(t,x) - f(t,x), \\
		\pt\left(\beta(t,x) v(t,x)\right)   & = - \curl u(t,x) - \tilde\rho(t,x)v(t,x) \hspace{1pt}- g(t,x),
	\end{split}
	\hspace{0.5cm}                                    &                                                                                       & t\geq 0, x \in \Omega_\pm,                                                 \\
	u(t,x) \times \nu(x)                              & = \bndInh(t,x),                                                                       &                            & t\geq 0, x \in \partial \Omega,               \\
	\llbracket u(t,x) \times \nu^\Gamma(x) \rrbracket & =\kappa(t,x),\quad\llbracket v(t,x) \times \nu^\Gamma(x) \rrbracket = \varkappa(t,x),
	                                                  &                                                                                       & t\geq 0, x \in \Gamma,                                                     \\
	u(0,x)                                            & = u_0(x), \hspace{20mm} v(0,x)= v_0(x),                                               &                            & x \in \Omega_\pm,\label{eq:maxwell-original2}
\end{alignat}
where $\llbracket f \rrbracket \coloneqq \tr_{\partial \upsilon} f - \tr_{\partial \omega} f$ denotes the jump across the interface. We may also include a (purely mathematical) magnetic conductivity $\tilde \rho$.
For energy estimates it is favorable to commute $\alpha$ and $\beta$ with $\pt$, which yields a symmetric, hyperbolic system.
In this setting the well-posedness theory has been established in \cite{schnaubelt2022}.
We thus rewrite the above system as
\begin{alignat}{-1}
	\begin{split}\label{eq:maxwell-lem}
		\alpha(t,x) \pt u(t,x) & = \hphantom{-}\curl v(t,x) - \sigma(t,x)u(t,x) - f(t,x), \\
		\beta(t,x) \pt v(t,x)  & = - \curl u(t,x) - \rho(t,x)v(t,x) \hspace{1pt}- g(t,x),
	\end{split}
	\quad                                             &                                                                                       & t\geq 0, x \in \Omega_\pm,                                                                               \\
	u(t,x) \times \nu(x)                              & = \bndInh(t,x),                                                                       &                                                & t\geq 0, x \in \partial \Omega, \label{eq:maxwell-lem3} \\
	\llbracket u(t,x) \times \nu^\Gamma(x) \rrbracket & =\kappa(t,x),\quad\llbracket v(t,x) \times \nu^\Gamma(x) \rrbracket = \varkappa(t,x),
	                                                  &                                                                                       & t\geq 0, x \in \Gamma, \label{eq:maxwell-lem2}                                                           \\
	u(0,x)                                            & = u_0(x), \hspace{20mm} v(0,x)= v_0(x),                                               &                                                & x \in \Omega_\pm, \label{eq:maxwell-lem4}
\end{alignat}
where
\begin{equation}
	\label{eq:transform-tilde-sigma}
	\sigma = \tilde \sigma + \pt\alpha \quad \text{and}\quad \rho = \tilde \rho + \pt\beta.
\end{equation}
In the following we will only use the formulation \eqref{eq:maxwell-lem}--\eqref{eq:maxwell-lem4} of the Maxwell system.

We employ the spaces
\begin{align*}
	\cC^{\ell}              & \coloneqq \bigcap_{j=0}^{\ell}C^{\ell - j}\left([0,T], H^{j}_\Gamma(\Omega,\Rsym)\right) ,                                   \\
	\cC^{\ell}_\eta         & \coloneqq \left\{A \in \cC^{\ell} \mid A(t) \geq \eta \text{ for all }t \in [0,T]\right\}, \qquad \text{for } \ell \in \N_0, \\
	\|A\|_{\cC^{\ell}_\eta} & = \|A\|_{\cC^{\ell}}\coloneqq \sup_{\substack{j=0,\dots\ell                                                                  \\t\in[0,T]}} \smallnorm{\pt^{\ell - j}A(t)}_{H^{j}_\Gamma},
\end{align*}
with
\[
	H_\Gamma^j(\Omega,\Rsym) \coloneqq \left\{f \in L^2(\Omega,\Rsym) \Bigm\vert f|_\omega \in H^j(\omega,\Rsym) \text{ and } f|_\upsilon \in H^j(\upsilon,\Rsym)\right\}.
\]
The coefficients are required to satisfy
\begin{align}\label{eq:maxwell-lem-coeff}
	\alpha,\beta & \in \cC^{\tilde m}_\eta \text{ for some } \eta >0 \quad \text{and}\quad \sigma, \rho \in \cC^{\tilde m},
\end{align}
with $\tilde m\coloneqq \max\{m, 3\}$, cf.\ Section~2 in \cite{schnaubelt2022}.
Furthermore, we assume that
\begin{equation}
	\label{eq:additional-reg}
	\pt\alpha,\pt \beta \in \cC^{\hat m}
\end{equation}
with $\hat m\coloneqq \max\{m, 2\}$.
The latter we use to control error terms when estimating the divergence of $\alpha u$ and $\beta v$ in the proof of the convergence of the approximation by means of a div-curl lemma, see \eqref{eq:proof:ibp}.
Note that for $m \leq 2$, i.e.\ $\tilde m =3$ and $\hat m =2$, condition \eqref{eq:maxwell-lem-coeff} already implies \eqref{eq:additional-reg}.
In Remark~\ref{rem:quasilinear} we discuss the application of our results to a time-derived, quasilinear system in the context of condition~\eqref{eq:additional-reg}.
We further assume
\begin{equation}\label{eq:maxwell-lem-IVs}
	f,g \in H^m([0,T]\times \Omega_\pm),\ \bndInh\in T^m(\partial\Omega),\ \kappa,\varkappa \in J^m(\Gamma), \ \text{and} \ u_0, v_0 \in H_\Gamma^m(\Omega),
\end{equation}
where the space of spatial traces
\begin{equation*}
	T^m(\partial\Omega) \coloneqq \tr_{t,\partial\Omega} \cG^{m+1}_\Gamma,\qquad J^m(\Gamma) \coloneqq\llbracket \cG^{m+1}_\Gamma \times \nu^\Gamma\rrbracket,
\end{equation*}
are subsets of $E^m(\Xi) \coloneqq \cap_{j=0}^m H^j([0,T], H^{m+1/2-j}(\Xi))$, for $\Xi \in\{\partial\Omega,\Gamma\}$, see \cite{schnaubelt2022} and \cite{spitz2018}, allowing us to extend $\bndInh,\kappa,\varkappa$ to the interior as functions $\bndInh^\circ,\kappa^\circ,\varkappa^\circ$ in
\[
	\cG^{m+1}_\Gamma \coloneqq \cG_\Gamma^{m+1}([0,T]\times \Omega)\coloneqq \bigcap_{j=0}^{m+1} C^j([0,T], H_\Gamma^{m+1-j}(\Omega)).
\]
These spaces naturally occur, when considering tangential derivatives of solutions to the Maxwell system with homogeneous boundary and interface conditions as in \cite{nutt2026-2}.

\begin{remark}\label{rem:dense}
	a) Note that the coefficients $\alpha, \beta \in \cC^{\tilde m}_\eta$ satisfy
	\begin{align*}
		\alpha|_{\Xi}, \beta|_{\Xi}, \sigma|_{\Xi}, \rho|_{\Xi} \in & \{A \in W^{1,\infty}([0,T] \times \Xi)^{3 \times 3} \colon                                                                       \\
		                                                            & \partial^\alpha A \in L^\infty([0,T], L^2(\Xi)) \text{ for all } \alpha \in \N_0^4 \text{ with } 1 \leq |\alpha| \leq \tilde m\}
	\end{align*}
	for $\Xi \in \{\upsilon, \omega\}$ and that $\alpha,\beta$ are uniformly positive definite.
	Hence, we may apply the local well-posedness result from Theorem~3.1 in \cite{schnaubelt2022} to our system~\eqref{eq:maxwell-lem} above.

	b) By extending the coefficients, using the usual construction of Theorem~5.19 in \cite{adams2003} in time, and a mollification argument, one can show that the space $C^\infty([0,T] \times \ol{\Omega_\pm}, \Rsym)$ is dense in $C^{\tilde m - j}([0,T], H^{j}_\Gamma(\Omega,\Rsym))$. Thus, we can approximate our coefficients smoothly, see \ref{def:approx-data}.
\end{remark}

In order to study functions at the boundary and interface we use the vector-valued normal trace $\tr_n \colon H(\div) \to H^{-1/2}(\Xi)^3$ as well as the tangential traces $\tr_t,\tr_\tau \colon H(\curl) \to H^{-1/2}(\Xi)^3$\\[-0.2cm]
\noindent
\begin{minipage}[t]{.48\textwidth}
	\vspace{0.2cm}
	\begin{alignat*}{-1}
		\tr_n u    & = \tr_{n, \Xi} u   &  & \coloneqq (u \cdot \nu^\Xi) \nu^\Xi,         \\
		\tr_t u    & = \tr_{t, \Xi} u   &  & \coloneqq u \times \nu^\Xi ,                 \\
		\tr_\tau u & = \tr_{\tau,\Xi} u &  & \coloneqq \nu^\Xi \times (u \times \nu^\Xi),
	\end{alignat*}
\end{minipage}
\begin{minipage}[t]{.51\textwidth}
	\vspace{-0.4cm}
	\centering
	\resizebox{\linewidth}{!}{%
		\tdplotsetmaincoords{73}{90}

		\begin{tikzpicture}[scale=2,tdplot_main_coords]

			\coordinate (O) at (0,0,0);

			\coordinate (N) at (0,0,0.85);       
			\coordinate (U) at (0.6,0.6,0.6); 

			\coordinate (Unorm) at (0,0,0.6);

			\coordinate (UxN) at (0.6,-0.6,0);

			\coordinate (NxUxN) at (0.6,0.6,0);

			\draw[->, thin] (O) -- (N) node[anchor=south]{$\nu^\Xi$};

			\draw[->, thick] (O) -- (U) node[anchor=west]{$u$};

			\draw[->, thick] (O) -- (Unorm) node[anchor= east]{$(u \cdot \nu^\Xi)\nu^\Xi$};

			\draw[->, thick] (O) -- (UxN) node[anchor=south]{$u \times \nu^\Xi$};

			\draw[->, thick] (O) -- (NxUxN) node[anchor= west]{$\nu^\Xi \times (u \times \nu^\Xi)$};

			\draw[dashed] (U) -- (Unorm);
			\draw[dashed] (U) -- (NxUxN);

			\draw (-45:0.3) arc (-45:45:0.3);
			\fill (0:0.15) circle[radius=0.02]; 

			\tdplotsetrotatedcoords{-135}{-90}{45}
			\draw[tdplot_rotated_coords] (-45:0.3) arc (-45:45:0.3);
			\fill[tdplot_rotated_coords] (0:0.15) circle[radius=0.02]; 

			\shade[ball color=gray!10, opacity=0.3] (0,0,0) circle (1.1);

		\end{tikzpicture}
	}
	\vspace{-0.2cm}
\end{minipage}
where $\Xi$ is $\partial \upsilon$, $\partial \omega$ or $\partial \Omega$ and $\nu=\nu^\Xi$ denotes the corresponding outer normal.
See Chapter~2.2.2 in \cite{ciarlet18} or Chapter~2.1 in \cite{cessenat1996} for mapping properties and further details on these traces.

Finally, we introduce the notation ``$\lesssim$'', whenever an inequality holds up to a constant not
depending on time or the explicit norm of the solution.

\subsection*{Compatibility Conditions}

Assume for a moment that we have a solution $(u,v) \in \cG^m_\Gamma([0,T])$ of \eqref{eq:maxwell-lem} satisfying the boundary and interface conditions
\begin{equation*}
	\tr_{t, \partial \Omega} u(t) = \bndInh, \quad \llbracket u(t) \times \nu^\Gamma \rrbracket = \kappa, \quad \llbracket v(t) \times \nu^\Gamma \rrbracket = \varkappa,
\end{equation*}
for all $t \geq 0$.
Since the trace $\tr_{t, \partial \Omega}$ and the interface jump $\llbracket \cdot \times \nu^\Gamma \rrbracket$ act only on spatial variables, they commute with $\pt^i$ for $i<m$ resulting in
\[
	\tr_{t, \partial \Omega} \pt^i u(t) = \pt^i \bndInh, \quad \llbracket \pt^i u(t) \times \nu^\Gamma \rrbracket = \pt^i \kappa, \quad \llbracket \pt^i v(t) \times \nu^\Gamma \rrbracket = \pt^i \varkappa.
\]
Now, the system~\eqref{eq:maxwell-lem} relates time derivatives to spatial derivatives yielding the conditions
\begin{align*}
	\tr_{t, \partial \Omega} \left(\pt^{i-1} \left(\inv\alpha( \curl v - \sigma u - f) \right)\right)   & = \pt^i \bndInh,   \\
	\llbracket \pt^{i-1} \left(\inv\alpha( \curl v - \sigma u - f) \right) \times \nu^\Gamma \rrbracket & = \pt^i \kappa,    \\
	\llbracket \pt^{i-1} \left(\inv\beta (- \curl u - \rho v - g)\right) \times \nu^\Gamma \rrbracket   & = \pt^i \varkappa,
\end{align*}
for $1\leq i<m$.
By applying \eqref{eq:maxwell-lem} again we can recursively translate these conditions on the time derivatives into conditions on the spatial derivatives of $(u,v)$.
In particular, they have to hold at $t=0$, i.e., the initial values $u_0, v_0 \in H_\Gamma^m(\Omega)$ of a solution $(u,v) \in \cG_\Gamma^m$ to \eqref{eq:maxwell-lem}--\eqref{eq:maxwell-lem4} must satisfy additional \textit{compatibility conditions}.

The compatibility conditions of order $m$ are then defined as follows, see Equation~(2.6) in \cite{schnaubelt2022}.
As usual the sum from $1$ to $0$ is empty.
\begin{definition}[Compatibility Condition of order $m$]\label{def:CC}
	Let $(u_{(0)}, v_{(0)}) \coloneqq (u_0, v_0) \in H^{m}_\Gamma(\Omega)$. For $1 \leq i \leq m$, we define $u_{(i)}, v_{(i)}$ by recursively differentiating the system at $t=0$, obtaining
		{\setlength{\abovedisplayskip}{4pt}
			\setlength{\belowdisplayskip}{4pt}
			\begin{alignat*}{2}
				u_{(i)} & \coloneqq \alpha(0)^{-1} \Bigg( &  & \curl v_{(i-1)} - \sigma(0) u_{(i-1)} - \pt^{i-1} f(0)                                                                          \\*
				        &                                 &  & - \sum_{\ell=1}^{i-1} \binom{i-1}{\ell} \left((\pt^\ell \alpha)(0) u_{(i-\ell)}+\pt^\ell\sigma(0) u_{(i-1-\ell)} \right)\Bigg), \\*
				v_{(i)} & \coloneqq \beta(0)^{-1} \Bigg(  &  & -\curl u_{(i-1)} - \rho(0) v_{(i-1)}- \pt^{i-1} g(0)                                                                            \\*
				        &                                 &  & - \sum_{\ell=1}^{i-1} \binom{i-1}{\ell} \left((\pt^\ell \beta)(0) v_{(i-\ell)} +\pt^\ell\rho(0)v_{(i-1-\ell)} \right) \Bigg).
			\end{alignat*}}
	Then, the compatibility conditions of order $m$ are given by
	\begin{equation}\label{eq:def-CC-traces}
		\tr_{t, \partial \Omega} u_{(i)} = \pt^i \bndInh, \quad \llbracket u_{(i)} \times \nu^\Gamma \rrbracket = \pt^i \kappa \quad \text{and} \quad \llbracket v_{(i)} \times \nu^\Gamma \rrbracket = \pt^i \varkappa \qquad \text{for }i< m,
	\end{equation}
	which we denote by $(CC_{\alpha,\beta,\sigma,\rho,f,g,\bndInh,\kappa,\varkappa,m}) =(CC_m)$.
	The compatibility conditions $(CC_0)$ of order zero do not require any conditions.
\end{definition}

Note that the compatibility conditions $(CC_m)$ are conditions on the time derivatives up to order $m-1$. Hence, the traces in \eqref{eq:def-CC-traces} are well-defined in $H^{1/2}$.
Theorem~3.1 and the comment after (1.7) in \cite{schnaubelt2022} then provide the following well-posedness result and a-priori estimate.
\begin{proposition}\label{prop:lwp-apriori}
	Let $\eta > 0$, $m,k \in \N_0$ with $\tilde m + 2 \leq k$ and $\Omega$ as above. Fix $T, r >0$.
	Take coefficients $(\alpha,\beta,\sigma,\rho)$ satisfying \eqref{eq:maxwell-lem-coeff} with $\norm{(\alpha,\beta,\sigma,\rho)}_{\cC^{\tilde m}}<r$, and data $(u_0,v_0)$, $(\bndInh,\kappa,\varkappa)$ and $(f, g)$ satisfying \eqref{eq:maxwell-lem-IVs}.
	Let $(u_0,v_0)$ fulfill the compatibility conditions $(CC_m)$.
	Then the system~\eqref{eq:maxwell-lem}--\eqref{eq:maxwell-lem4} has a unique solution $(u,v) \in \cG_\Gamma^m([0,T])$, such that
	\begin{equation}
		\label{eq:lwp-apriori}
		\begin{aligned}
			\norm{(u,v)}_{\cG^m_\Gamma}^2 & \leq C_{r,T} \Bigl(\norm{(u_0,v_0)}_{H_\Gamma^{m}}^2 +\norm{(f,g)}_{H^m([0,T]\times \Omega_\pm)}^2                        \\
			                              & \seq \hphantom{C_{r,T} \Bigl(}+ \norm{\bndInh}_{E^m(\partial\Omega)}^2 + \norm{(\kappa,\varkappa)}_{E^m(\Gamma)}^2\Bigr).
		\end{aligned}
	\end{equation}
	for a constant $C_{r,T}>0$ depending on $r$ and $T$, but not on $u_0,v_0,f,g$ or $\bndInh,\kappa,\varkappa$.
\end{proposition}

\section{Approximation of Data}

We first show that we can always approximate first-order compatibility conditions. These are much simpler than higher order compatibility conditions, since they do not involve derivatives of other fields.
\begin{lemma}\label{lem:approx-CC1}
	Let $m,a \in \N_0$ with $m+a \leq k$, $\bndInh \in T^{m+a}(\partial\Omega)$, and $\kappa,\varkappa \in J^{m+a}(\Gamma)$. Then the set
	\[
		H_{\Gamma, CC_1}^{m+a}(\Omega) \coloneqq \left\{(u,v)\in H_\Gamma^{m+a}(\Omega) \bigm\vert \tr_{t,\partial \Omega} u = \bndInh,\ \left\llbracket (u \times \nu^\Gamma, v \times \nu^\Gamma) \right\rrbracket = (\kappa,\varkappa)\right\}
	\]
	is dense in $L^2(\Omega)$ as well as in
	\[
		H_{\Gamma, CC_1}^{m}(\Omega) = \left\{(u,v)\in H_\Gamma^m(\Omega) \bigm\vert \tr_{t,\partial \Omega} u = \bndInh,\ \left\llbracket (u \times \nu^\Gamma, v \times \nu^\Gamma) \right\rrbracket = (\kappa,\varkappa)\right\}.
	\]
\end{lemma}

This means that we can approximate adding or preserving the compatibility condition of first-order.

\begin{proof}
	We first extend $(\bndInh,\kappa)$ and $(0,\varkappa)$ to functions $W^u,W^v \in H^{m+a}(\omega)$ satisfying
	\[
		\tr_{t,\Gamma}(W^u,W^v) = (\kappa,\varkappa) \quad \text{and}\quad \tr_{t,\partial\Omega}(W^u,W^v) = (\bndInh,0)
	\]
	simultaneously by using cutoffs $\chi_{\Xi} \in C_c^\infty(\bar\omega)$, for $\Xi \in\{\Gamma, \partial\Omega\}$, such that $\chi_{\Xi}|_{\Xi} = 1$ and $\supp \chi_{\Gamma} \cap \supp \chi_{\partial \Omega} = \emptyset$ and setting
	\begin{equation}
		\label{eq:simultaneous-extensions}
		W^u = \chi_{\partial\Omega} \bndInh^\circ + \chi_{\Gamma} \kappa^\circ \quad \text{and} \quad W^v = \chi_{\Gamma} \varkappa^\circ.
	\end{equation}
	We then extend $W^u$ and $W^v$ by $0$ to $\upsilon$ and denote the extensions also by $W^u,W^v \in H_\Gamma^{m+a}(\Omega_\pm)$.
	The density in $L^2$ is then clear since $(W^u,W^v) + C_c^\infty(\Omega_\pm)$ is dense in $L^2(\Omega_\pm) = L^2(\Omega)$.

	Next, we want to approximate $(u,v) \in H_{\Gamma, CC_1}^{m}(\Omega)$ by $(u_n,v_n) \in H^{m+a}_\Gamma(\Omega)$ satisfying the compatibility condition of order~$1$,
	namely
	\begin{equation}\label{eq:proof-comp-cond-1st-order}
		\tr_{t,\partial \Omega}(u_n) = \bndInh, \qquad \left\llbracket (u_n \times \nu^\Gamma, v_n \times \nu^\Gamma) \right\rrbracket = (\kappa,\varkappa).
	\end{equation}

	For the construction of these maps, we choose $(u^\upsilon_n,v^\upsilon_n) \in H^{m+a}(\upsilon)$ such that $(u^\upsilon_n,v^\upsilon_n) \to (u,v)$ in $H^m(\upsilon)$ and set
	\[
		(u_n,v_n)|_\upsilon \coloneqq (u^\upsilon_n,v^\upsilon_n).
	\]
	We next approximate $(u,v)$ in $H^m(\omega)$ by $(u^\omega_n,v^\omega_n) \in H^{m+a}(\omega)$.
	Now we have to adjust the trace terms in order to satisfy the compatibility conditions.
	The map
	\[
		c_n \coloneqq \begin{cases}
			\left(\tr_{t,\partial \Omega}(u^\omega_n) - \bndInh, 0\right), \quad                                                                 & \text{on } \partial \Omega,   \\
			\tr_{t,\Gamma}\left((u^\omega_n,v^\omega_n)\right)-\tr_{t,\Gamma}\left((u^\upsilon_n,v^\upsilon_n)\right) + (\kappa,\varkappa),\quad & \text{on } \partial \upsilon,
		\end{cases}
	\]
	then belongs to $\big(H^{m+a-1/2}(\partial \omega)\big)^6$.
	Applying the extension operator $R$ componentwise, we then define
	\[
		(u_n,v_n)|_\omega \coloneqq (u^\omega_n,v^\omega_n) - R c_n.
	\]
	Clearly, $(u_n,v_n) \in H^{m+a}_\Gamma$ satisfies the compatibility condition~\eqref{eq:proof-comp-cond-1st-order}. Since $(u,v)$ satisfies the compatibility conditions by definition, we have $\tr_{t,\partial \Omega} u = \bndInh$ and $\tr_{t,\partial \upsilon}(u,v) - \tr_{t,\partial \omega}(u,v)= (\kappa,\varkappa)$.
	Substituting these into our bound for $c_n$ allows us to split the jump difference back into the separate domain errors. So, we can estimate
	\begin{align*}
		 & \norm{(u,v)-(u_n,v_n)}_{H_\Gamma^m}                                                                                                                                                                                                      \\
		 & \leq \norm{(u,v)-(u^\upsilon_n,v^\upsilon_n)}_{H^m(\upsilon)} + \norm{(u,v)-(u^\omega_n,v^\omega_n)}_{H^m(\omega)} + \norm{R c_n}_{H^m(\omega)}                                                                                          \\
		 & \lesssim \norm{(u,v)-(u^\upsilon_n,v^\upsilon_n)}_{H^m(\upsilon)} + \norm{(u,v)-(u^\omega_n,v^\omega_n)}_{H^m(\omega)}                                                                                                                   \\*
		 & \seq + \norm{\tr_{t,\partial \Omega}\left(u^\omega_n\right) - \bndInh}_{H^{m-1/2}(\partial \Omega)}                                                                                                                                      \\*
		 & \seq + \norm{\tr_{t,\partial \omega}\left(u^\omega_n,v^\omega_n\right) - \tr_{t,\partial \upsilon}\left(u^\upsilon_n,v^\upsilon_n\right) + (\kappa,\varkappa)}_{H^{m-1/2}(\partial \upsilon)}                                            \\
		 & = \norm{(u,v)-(u^\upsilon_n,v^\upsilon_n)}_{H^m(\upsilon)} + \norm{(u,v)-(u^\omega_n,v^\omega_n)}_{H^m(\omega)}                                                                                                                          \\
		 & \seq + \norm{\tr_{t,\partial \Omega}\left(u^\omega_n - u\right)}_{H^{m-1/2}(\partial \Omega)}                                                                                                                                            \\
		 & \seq + \norm{\tr_{t,\partial \omega}\left(u^\omega_n,v^\omega_n\right) - \tr_{t,\partial \omega}(u,v) + \tr_{t,\partial \upsilon}(u,v) - \tr_{t,\partial \upsilon}\left(u^\upsilon_n,v^\upsilon_n\right)}_{H^{m-1/2}(\partial \upsilon)} \\
		 & \lesssim \norm{(u,v)-(u^\upsilon_n,v^\upsilon_n)}_{H^m(\upsilon)} + \norm{(u,v)-(u^\omega_n,v^\omega_n)}_{H^m(\omega)}                                                                                                                   \\
		 & \to 0 \quad \text{as } n\to\infty. \qedhere
	\end{align*}
\end{proof}

\begin{remark}\label{rem:approx-m-geq-1}
	For a given solution $(u,v)$ of \eqref{eq:maxwell-lem}--\eqref{eq:maxwell-lem4} in $\cG_\Gamma^0([0,T])$, sequences $(f_n,g_n) \in H^{1}([0,T] \times \Omega_\pm)$ such that $(f_n, g_n)\to (f,g)$ in $L^2([0,T]\times\Omega)$ and $\bndInh_n \in T^1(\partial\Omega), (\kappa_n,\varkappa_n) \in J^1(\Gamma)$ such that $\bndInh_n \to \bndInh$ in $T^0(\partial\Omega)$ and $(\kappa_n,\varkappa_n) \to (\kappa,\varkappa)$ in $J^0(\Gamma)$, Lemma~\ref{lem:approx-CC1} allows us to approximate its initial values $(u_0,v_0) \in L^2(\Omega)$ by $(u_{0,n_k}, v_{0,n_k})$ in $H_{\Gamma, CC_1}^{1}(\Omega)$ with respect to $(\bndInh_n,\kappa_n,\varkappa_n)$.

	By the a-priori estimate \eqref{eq:lwp-apriori}, the corresponding diagonal sequence of solutions $(u_{n_n},v_{n_n}) \in \cG_\Gamma^1([0,T])$ of \eqref{eq:maxwell-lem}--\eqref{eq:maxwell-lem4} converge to $(u,v)$ in $\cG_\Gamma^0([0,T])$.
\end{remark}

In view of Remark~\ref{rem:approx-m-geq-1}, we assume that $m \geq 1$ in the following.
In order to obtain solutions of higher regularity than $\tilde m$, we need additional regularity the coefficients and inhomogeneities, so that the coefficients also have to be approximated. These functions exists by Remark~\ref{rem:dense}. They also enter the compatibility conditions of the approximation.
\medskip
\begin{enumerate}[leftmargin=30pt, label=\textup{($\text{App}$)}, after=\bigskip]
	\item \label{def:approx-data}\quad \parbox[m]{\linegoal}{
		      We choose $(\alpha_n, \beta_n, \sigma_n, \rho_n) \in C^\infty\left([0,T] \times \ol{\Omega_\pm}; \Rsym\right)^4 \cap \left(\cC^{\tilde{m}}_\eta\times\cC^{\tilde{m}}_\eta\times\cC^{\tilde{m}}\times\cC^{\tilde{m}}\right) $ with $\alpha_n \to \alpha$, $\beta_n \to \beta$ in $\cC^{\tilde{m}}$, $\pt\alpha_n \to \pt\alpha$, $\pt\beta_n \to \pt\beta$ in $\cC^{\hat{m}}$, $\sigma_n \to \sigma$, and $\rho_n \to \rho$
		      in $\cC^{\tilde{m}}$, $(f_n,g_n) \in C^\infty([0,T] \times \overline{\Omega_\pm})$ with $f_n \to f$, $g_n \to g$ in $H^{m}([0,T] \times \Omega_\pm)$.
		      Furthermore, we approximate the $\cG^{m+1}_\Gamma$-extensions $\bndInh^\circ, \kappa^\circ$ and $\varkappa^\circ$ by smooth functions $\bndInh_n^\circ, \kappa_n^\circ, \varkappa_n^\circ \in C^\infty([0,T] \times \overline{\Omega_\pm})$ in $\cG^{m+1}_\Gamma$.
		      We set $\bndInh_n = \tr_{\partial\Omega} \bndInh_n^\circ\in C^\infty([0,T] \times \partial\Omega)$. Then $\bndInh_n \to \bndInh$ in $E^m(\partial\Omega)$. Similarly, we set $\kappa_n \coloneqq \llbracket \kappa_n^\circ \rrbracket, \varkappa_n \coloneqq \llbracket\varkappa_n^\circ \rrbracket \in C^\infty([0,T] \times \Gamma)$ such that $(\kappa_n,\varkappa_n)\to (\kappa,\varkappa)$ in $E^m(\Gamma)$.}
\end{enumerate}

In a core step, we start from fields $(u_0,v_0)$ which satisfy the $m$-th order compatibility condition up to an error of size $r$.
The next lemma then provides fields fulfilling the $(m+a)$-th order compatibility condition whose distance to $(u_0,v_0)$ in $H^m_\Gamma$ is less than $Cr$.
\begin{proposition}\label{prop:fix-CC}
	Let $a \in \N_0$, $m \geq 1$ with $m+a\leq k$ and $r>0$. Let $(u^{\text{app}}_{(0)},v^{\text{app}}_{(0)}) \coloneqq (u_0,v_0)\in H^{m+a}_{\Gamma, CC_1}(\Omega)$, and set recursively
		{\setlength{\abovedisplayskip}{4pt}
			\setlength{\belowdisplayskip}{4pt}
			\begin{alignat*}{2}
				u^{\text{app}}_{(i)} & = \alpha_n(0)^{-1} \Bigg( &  & \curl v^{\text{app}}_{(i-1)} - \sigma_n(0) u^{\text{app}}_{(i-1)} - \pt^{i-1} f_n(0)                                                                               \\
				                     &                           &  & - \sum_{\ell=1}^{i-1} \binom{i-1}{\ell} \left((\pt^\ell \alpha_n)(0) u^{\text{app}}_{(i-\ell)} + (\pt^\ell \sigma_n)(0) u^{\text{app}}_{(i-1-\ell)}\right) \Bigg), \\
				v^{\text{app}}_{(i)} & = \beta_n(0)^{-1} \Bigg(  &  & -\curl u^{\text{app}}_{(i-1)} - \rho_n(0) v^{\text{app}}_{(i-1)} - \pt^{i-1} g_n(0)                                                                                \\
				                     &                           &  & - \sum_{\ell=1}^{i-1} \binom{i-1}{\ell} \left((\pt^\ell \beta_n)(0) v^{\text{app}}_{(i-\ell)} + (\pt^\ell \rho_n)(0) v^{\text{app}}_{(i-1-\ell)}\right)\Bigg),
			\end{alignat*}
		}
	for $1\leq i \leq m$.
	Assume that	$(u_0,v_0)$ satisfy the compatibility conditions of order $m$ approximately, i.e.,
	\begin{align}
		\label{eq:approx-comp-cond}
		\notag\norm{\tr_{t, \partial \Omega} u^{\text{app}}_{(i)} - \pt^i \bndInh_n(0)}_{H^{m-i-1/2}(\partial\Omega)}           & < r, \\*
		\norm{\llbracket u^{\text{app}}_{(i)} \times \nu^\Gamma \rrbracket -\pt^i \kappa_n(0)}_{H^{m-i-1/2}(\Gamma)}            & < r, \\*
		\notag \norm{\llbracket v^{\text{app}}_{(i)} \times \nu^\Gamma \rrbracket - \pt^i \varkappa_n(0)}_{H^{m-i-1/2}(\Gamma)} & < r,
	\end{align}
	for $0\leq i <m$.
	Then, there exists a constant $C>0$, independent of $r$, and functions $(u^{m+a}_0, v^{m+a}_0)\in H^{m+a}_\Gamma(\Omega)$ satisfying $(CC_{\alpha_n, \beta_n, \sigma_n, \rho_n, f_n, g_n, \bndInh_n, \kappa_n, \varkappa_n, m+a})$ such that
	\begin{equation}\label{eq:approx-close}
		\norm{(u_0,v_0) - (u^{m+a}_0, v^{m+a}_0)}_{H^{m}_\Gamma} \leq C r.
	\end{equation}
\end{proposition}

For a given $r>0$ the condition \eqref{eq:approx-comp-cond} can be satisfied when approximating coefficients and initial data satisfying the compatibility condition of order $m$ as in \ref{def:approx-data} for large $n$ and then applying Lemma~\ref{lem:approx-CC1}.

\begin{proof}[Proof of Proposition~\ref{prop:fix-CC}]
	We proceed iteratively over the order $j$ of the compatibility condition from $0$ to $m$. The cases $m+\hat{\jmath}$ (for $0 < \hat{\jmath}\leq a$) require a different argument at the end.
	\begin{enumerate}[label = \arabic*), wide, labelwidth=0pt, labelindent=0pt, itemsep=0.5em, topsep=0pt]
		\item Assume that for some $0 < j<m$, the pair $(u_0^j, v_0^j) \in H^{m+a}_{\Gamma}(\Omega)$ satisfies the exact compatibility condition $(CC_j)$ and the estimate $\smallnorm{(u_0,v_0) - (u_0^j, v_0^j)}_{H^{m}(\Omega)} \leq C r$.
		      We construct $(u_0^{j+1}, v_0^{j+1})$ that fulfill the condition of order $j+1$, while preserving the conditions of order $k$ for $k < j$ by adding a correction term $(c_0^j, d_0^j) \in H^{m+a}(\omega) \cap H^j_0(\omega)$ (and $0$ on $\upsilon$) to $(u_0^j, v_0^j)$.
		      Set $(u_0^{j+1}, v_0^{j+1}) \coloneqq (u_0^j, v_0^j) + (c_0^j, d_0^j)$.

		      We first want to study the effect of adding $(c_0^j, d_0^j) \in H^{m+a}(\omega) \cap H^j_0(\omega)$. We will choose it explicitly later on.
		      For $i<j$ the compatibility conditions
		      \begin{equation*}
			      \tr_{t, \partial \Omega} u^{j+1}_{(i)} = \pt^i \bndInh_n(0), \quad \left\llbracket u^{j+1}_{(i)} \times \nu^\Gamma \right\rrbracket = \pt^i \kappa_n(0) \quad \text{and} \quad \left\llbracket v^{j+1}_{(i)} \times \nu^\Gamma \right\rrbracket = \pt^i \varkappa_n(0)
		      \end{equation*}
		      only depend on the traces of derivatives of order $\leq i$ of $(u_0^{j+1}, v_0^{j+1})$.
		      The fields $(c_0^j, d_0^j)\in H^j_0(\omega)$ do not influence these traces for $i<j$. So, we obtain $\tr_{t, \Xi}(u_{(i)}^{j+1}, v_{(i)}^{j+1}) = \tr_{t, \Xi}(u_{(i)}^{j}, v_{(i)}^{j})$ for $\Xi \in \{\partial \Omega, \partial \omega \cap \Omega, \partial \upsilon\}$ and $i<j$.

		      The $j$-th compatibility condition of $(CC_{j+1})$ for $u^{j+1}_{(j)}$ imposes on $\curl v^{j+1}_{(j-1)}$ the identity
		      \begin{align*}
			       & \left\llbracket u^{j+1}_{(j)} \times \nu^\Gamma \right\rrbracket \stackrel{!}{=} \pt^j \kappa_n(0) \iff \left\llbracket \big(\alpha_n(0)^{-1} \curl v^{j+1}_{(j-1)} \big)\times \nu^\Gamma\right\rrbracket \stackrel{!}{=} I^v_j \\*
			       &
			      \begin{aligned}
				      \text{with } I^v_j = & I^v_j\begin{pmatrix}u_0^{j+1}\\v_0^{j+1}\end{pmatrix}\coloneqq \pt^j \kappa_n(0) + \Bigg\llbracket \alpha_n(0)^{-1} \Bigg( \sigma_n(0) u^{j+1}_{(j-1)} + \pt^{j-1} f_n(0)                   \\
				                           & \qquad+ \sum_{\ell=1}^{j-1} \binom{j-1}{\ell}\left( (\pt^\ell \alpha_n)(0) u^{j+1}_{(j-\ell)} +(\pt^\ell \sigma_n)(0) u^{j+1}_{(j-1-\ell)} \right) \Bigg) \times \nu^\Gamma\Bigg\rrbracket.
			      \end{aligned}
		      \end{align*}
		      Likewise, the compatibility condition for $v^{j+1}_{(j)}$ requires the equation
		      \begin{align*}
			       & \left\llbracket v^{j+1}_{(j)} \times \nu^\Gamma \right\rrbracket \stackrel{!}{=} \pt^j \varkappa_n(0) \iff \left\llbracket\big(\beta_n(0)^{-1} \curl u^{j+1}_{(j-1)}\big) \times \nu^\Gamma \right\rrbracket \stackrel{!}{=} I^u_j \\
			       & \begin{aligned}\text{with } I^u_j = & I^u_j\begin{pmatrix}u_0^{j+1}\\v_0^{j+1}\end{pmatrix} \coloneqq \pt^j \varkappa_n(0) -\Bigg\llbracket \beta_n(0)^{-1} \Bigg( \rho_n(0) v^{j+1}_{(j-1)} + \pt^{j-1} g_n(0) \\
                                    & \qquad+ \sum_{\ell=1}^{j-1} \binom{j-1}{\ell}\left( (\pt^\ell \beta_n)(0) v^{j+1}_{(j-\ell)} +(\pt^\ell \rho_n)(0) v^{j+1}_{(j-1-\ell)} \right) \Bigg) \times \nu^\Gamma\Bigg\rrbracket.
			         \end{aligned}
		      \end{align*}
		      At the outer boundary we obtain
		      \begin{align*}
			       & \tr_{t, \partial \Omega} u^{j+1}_{(j)} \stackrel{!}{=} \pt^j \bndInh_n(0) \iff \tr_{t, \partial \Omega} \left(\alpha_n(0)^{-1} \curl v^{j+1}_{(j-1)}\right) \stackrel{!}{=} B^v_j \\*
			       &
			      \begin{aligned}
				      \text{with } B^v_j = B^v_j\begin{pmatrix}u_0^{j+1}\\v_0^{j+1}\end{pmatrix} & \coloneqq \pt^j \bndInh_n(0) + \tr_{t, \partial \Omega}\Bigg( \alpha_n(0)^{-1} \Bigg( \sigma_n(0) u^{j+1}_{(j-1)} + \pt^{j-1} f_n(0)                        \\
				                                                                                 & + \sum_{\ell=1}^{j-1} \binom{j-1}{\ell}\left( (\pt^\ell \alpha_n)(0) u^{j+1}_{(j-\ell)} +(\pt^\ell \sigma_n)(0) u^{j+1}_{(j-1-\ell)} \right) \Bigg)\Bigg) .
			      \end{aligned}
		      \end{align*}
		      Note that the right-hand sides $B^v_j, I^u_j$, and $I^v_j$ do not depend on $(c_0^j,d_0^j)$, since the traces of $(u^{j+1}_{(i)},v^{j+1}_{(i)})$ for $i<j$ appearing therein remain unaffected by the correction.
		      So, to summarize we have to choose $(c_0^j,d_0^j)\in H^j_0(\omega)$ such that
		      \begin{align}\label{eq:proof:CC-for-u(j)}
			      \notag\tr_{t, \partial \Omega} \left(\alpha_n(0)^{-1} \curl v^{j+1}_{(j-1)}\right)   & \stackrel{!}{=} B^v_j(u_0^{j}, v_0^{j}), \\
			      \llbracket (\alpha_n(0)^{-1} \curl v^{j+1}_{(j-1)} )\times \nu^\Gamma\rrbracket      & \stackrel{!}{=} I^v_j(u_0^{j}, v_0^{j}), \\
			      \notag\llbracket(\beta_n(0)^{-1} \curl u^{j+1}_{(j-1)}) \times \nu^\Gamma \rrbracket & \stackrel{!}{=} I^u_j(u_0^{j}, v_0^{j}).
		      \end{align}

		      Let us take a closer look at the dependence of the left-hand side on $(c_0^j,d_0^j)$.
		      Recalling the iterative Definition~\ref{def:CC} of $(u^{j+1}_{(j-1)},v^{j+1}_{(j-1)})$, we realize that
		      \begin{equation}\label{eq:proof:split-principle-remainder}
			      \begin{pmatrix}
				      \alpha_n(0)^{-1} \curl v^{j+1}_{(j-1)} \\
				      \beta_n(0)^{-1} \curl u^{j+1}_{(j-1)}
			      \end{pmatrix} =
			      P^j\begin{pmatrix}u_0^{j+1}\\v_0^{j+1}\end{pmatrix} + R^j\begin{pmatrix}u_0^{j+1}\\v_0^{j+1}\end{pmatrix}
		      \end{equation}
		      where the principal part is given for even $j$ by
		      \begin{align*}
			      P^j(u,v) & = (-1)^{j/2} \begin{pmatrix}
				                              \bigl(\Pi_{i=1}^{j/2} \alpha_n^{-1}(0) \curl \beta_n^{-1}(0) \curl\bigr)u \\
				                              -\bigl(\Pi_{i=1}^{j/2} \beta_n^{-1}(0) \curl \alpha_n^{-1}(0) \curl\bigr) v
			                              \end{pmatrix}
			      \intertext{and for odd $j$}
			      P^j(u,v) & = (-1)^{(j-1)/2}\begin{pmatrix}
				                                 \bigl(\Pi_{i=1}^{(j-1)/2} \alpha_n^{-1}(0) \curl \beta_n^{-1}(0)\curl \bigr) \alpha_n^{-1}(0)\curl v \\
				                                 \bigl(\Pi_{i=1}^{(j-1)/2} \beta_n^{-1}(0) \curl \alpha_n^{-1}(0) \curl\bigr) \beta_n^{-1}(0) \curl u
			                                 \end{pmatrix}.
		      \end{align*}
		      The point being that the principal part $P^j$ is just a composition of curls and matrices, whereas the remainder $R^j (u_0^{j+1}, v_0^{j+1}) = R^j (u_0^{j+1}, v_0^{j+1}, \dots, D^{j-1}u_0^{j+1}, D^{j-1}v_0^{j+1})$ only depends in lower order on $\begin{pmatrix}u_0^{j+1}\\v_0^{j+1}\end{pmatrix}$. (Note that the remainder also depends on $(f_n,g_n)$.)
		      Since $(c^j_0,d^j_0) \in H^{j}_0(\Omega)$, the trace thus satisfies
		      \[
			      \tr_{t, \Xi} R^j(u_0^{j+1}, v_0^{j+1}) =\tr_{t, \Xi} R^j(u_0^j + c_0^j, v_0^j+ d_0^j)=\tr_{t, \Xi} R^j(u_0^j, v_0^j)
		      \]
		      for $\Xi \in \{\partial \Omega, \partial \omega, \partial \upsilon\}$.
		      Therefore, the linearity of $P^j$ and \eqref{eq:proof:split-principle-remainder} for $j$ imply
		      \begin{align*}
			      \tr_{t, \Xi} \begin{pmatrix}
				                   \alpha_n(0)^{-1} \curl v^{j+1}_{(j-1)} \\
				                   \beta_n(0)^{-1} \curl u^{j+1}_{(j-1)}
			                   \end{pmatrix} & =
			      \tr_{t, \Xi} \left(P^j(u^{j+1}_{0}, v^{j+1}_{0}) + R^j(u^{j}_{0}, v^{j}_{0})\right)                                                                    \\
			                                                & = \tr_{t, \Xi} \left(P^j(c^{j}_{0}, d^{j}_{0}) \right) + \tr_{t, \Xi} \begin{pmatrix}
				                                                                                                                        \alpha_n(0)^{-1} \curl v^{j}_{(j-1)} \\
				                                                                                                                        \beta_n(0)^{-1} \curl u^{j}_{(j-1)}
			                                                                                                                        \end{pmatrix}.
		      \end{align*}
		      Thus, the compatibility condition \eqref{eq:proof:CC-for-u(j)} becomes
		      \begin{alignat}{-1}\label{eq:proof:CC-for-u(j)-2}
			      \tr_{t, \partial \Omega} \left(P^j(c^{j}_{0}, d^{j}_{0}) \right)_u & \stackrel{!}{=} B^v_j(u_0^{j}, v_0^{j}) - \tr_{t, \partial \Omega} \left(\alpha_n(0)^{-1} \curl v^{j}_{(j-1)}\right)\notag                      \\
			                                                                         & = \pt^j \bndInh_n(0) - \tr_{t, \partial \Omega} u^j_{(j)}, \notag                                                                               \\
			      \llbracket P^j(c^{j}_{0}, d^{j}_{0})\times \nu^\Gamma \rrbracket   & \stackrel{!}{=} \begin{pmatrix}
				                                                                                           I^v_j(u_0^{j}, v_0^{j}) - \llbracket (\alpha_n(0)^{-1} \curl v^{j}_{(j-1)} )\times \nu^\Gamma\rrbracket \\
				                                                                                           I^u_j(u_0^{j}, v_0^{j}) - \llbracket(\beta_n(0)^{-1} \curl u^{j}_{(j-1)}) \times \nu^\Gamma \rrbracket
			                                                                                           \end{pmatrix} \notag \\
			                                                                         & =
			      \begin{pmatrix}
				      \pt^j \kappa_n(0) - \llbracket u^j_{(j)} \times \nu^\Gamma \rrbracket \\
				      \pt^j \varkappa_n(0) - \llbracket v^j_{(j)} \times \nu^\Gamma \rrbracket
			      \end{pmatrix},
		      \end{alignat}
		      where $(\cdot)_u$ denotes the first three components corresponding to the compatibility condition on $\tr_{t, \partial \Omega} u^{j+1}_{(j)}$.

		      Moreover, by Lemma~\ref{lem:Appendix:tr-curl} below, the (full) trace of the principal part is given by
		      \[
			      \tr_\Xi \left(P^j \begin{pmatrix}
					      c^{j}_{0} \\ d^{j}_{0}
				      \end{pmatrix}\right) = (-1)^{\lfloor{j/2}\rfloor}\begin{pmatrix}
				      \Pi^j_{i=1} T_{\zeta_i, x} & 0 \\ 0 & \Pi^j_{i=1} T_{\xi_i, x}
			      \end{pmatrix} \begin{pmatrix}
				      1 - \gamma & \gamma \\ -\gamma & 1 -\gamma
			      \end{pmatrix} \begin{pmatrix}
				      \tr(\partial_\nu^j c_0^j) \\ \tr(\partial_\nu^j d_0^j)
			      \end{pmatrix},
		      \]
		      where $\gamma = j \mod 2$, $T_{\zeta_i, x},T_{\xi_i, x}$ are given by \eqref{eq:def-T_alpha,x}, and $\zeta_i, \xi_i \in \{\alpha_n^{-1}(0) , \beta_n^{-1}(0)\}$ are smooth, uniformly positive definite matrices.
		      We want to solve for $\begin{pmatrix}
				      \tr(\partial_\nu^j c_0^j) \\ \tr(\partial_\nu^j d_0^j)
			      \end{pmatrix}$ after taking the tangential trace in the equality above. For this, note that the linear map
		      \[
			      A_{\tau}(x) \begin{pmatrix}
				      u \\v
			      \end{pmatrix} \coloneqq (-1)^{\lfloor{j/2}\rfloor}\left(
			      \begin{psmallmatrix}
					      \tr_{t, \Xi}& 0 \\ 0& \tr_{t, \Xi}
				      \end{psmallmatrix}
			      \begin{psmallmatrix}
					      \Pi^j_{i=1} T_{\zeta_i, x} & 0 \\ 0 & \Pi^j_{i=1} T_{\xi_i, x}
				      \end{psmallmatrix}	\begin{psmallmatrix}
					      1 - \gamma & \gamma \\ -\gamma & 1 -\gamma
				      \end{psmallmatrix}\begin{pmatrix}
					      u \\v
				      \end{pmatrix} \right)
		      \]
		      is invertible on the tangential plane for each $x \in \Xi$ for $\Xi \in \{\partial \Omega, \partial \omega \cap \Omega\}$ since the matrices $\zeta_i$ and $\xi_i$ are positive definite.
		      Let us denote the inverse on the tangential plane as
		      \[
			      A_{\tau}^{-1}(x)\colon \left(\langle\{\nu(x)\}\rangle^\perp\right)^2 \to \left(\langle\{\nu(x)\}\rangle^\perp\right)^2\subseteq\R^6.
		      \]
		      Furthermore, the matrices $\zeta_i$ and $\xi_i$ are smooth and the normal lies in $C^k(\partial \omega)$ because of $\partial\omega \in C^{k+1}$.
		      Thus, $A_{\tau}^{-1}$ boundedly maps $\{w \in H^\ell(\partial \omega) \mid w \cdot \nu = w \cdot \nu^\Gamma = 0\}$ into itself for $\ell \in [1, k]$.

		      Hence, when setting $(c_0^j, d_0^j)|_\upsilon = 0$ the compatibility condition in \eqref{eq:proof:CC-for-u(j)-2} is satisfied if and only if
		      \begin{align}\label{eq:proof:CC-for-u(j)-3}
			      \begin{aligned}
				      \tr_{\partial\Omega}\begin{pmatrix}
					                          \partial_\nu^j c_0^j \\ \partial_\nu^j d_0^j
				                          \end{pmatrix}    & \stackrel{!}{=} A_{\tau}^{-1}\begin{pmatrix} \pt^j \bndInh_n(0) - \tr_{t,\partial \Omega} u^j_{(j)}\\ 0 \end{pmatrix}, \\
				      \tr_{\partial\omega \cap \Omega}\begin{pmatrix}
					                                      \partial_\nu^j c_0^j \\ \partial_\nu^j d_0^j
				                                      \end{pmatrix} & \stackrel{!}{=} - A_{\tau}^{-1}\begin{pmatrix}
					                                                                                     \pt^j \kappa_n(0) - \llbracket u^j_{(j)} \times \nu^\Gamma \rrbracket \\
					                                                                                     \pt^j \varkappa_n(0) - \llbracket v^j_{(j)} \times \nu^\Gamma \rrbracket
				                                                                                     \end{pmatrix},
			      \end{aligned}
		      \end{align}
		      where the minus sign in the second equality stems from the definition of $\llbracket\cdot\rrbracket= \tr_{\partial \upsilon} - \tr_{\partial \omega}$.

		      As $\partial \omega \in C^{k+1}$, by Theorem~1.5.1.2 in \cite{grisvard1985} the map
		      \begin{equation}\label{eq:normal-trace-invertible}
			      \begin{aligned}
				      T\colon H^{m+a}(\omega) & \to \Pi_{s=0}^{j}H^{m+a - s - 1/2}(\partial \omega)\,;               \\
				      u                       & \mapsto (\tr u , \tr(\partial_\nu u) , \dots, \tr(\partial_\nu^j u))
			      \end{aligned}
		      \end{equation}
		      is surjective with continuous right inverse $T_{-1}$. Such an inverse can be explicitly constructed, see for example the proof of Lemma~3.36 and Theorem~3.37 in \cite{mclean2000}. One can then show that
		      \[
			      \norm{T_{-1} (\tr u , \tr(\partial_\nu u) , \dots, \tr(\partial_\nu^j u))}_{H^{m+a - \ell}(\omega)} \lesssim \sum_{s=0}^{j} \norm{\tr(\partial_\nu^{s} u)}_{H^{m+a - s - 1/2 - \ell}(\partial \omega)},
		      \]
		      for $\ell \in\{0, \dots, m+a-j-1\}$.
		      Hence, we find $c^j_0, d^j_0 \in H^{m+a}(\omega)\cap H^j_0(\omega)$ satisfying \eqref{eq:proof:CC-for-u(j)-3}.
		      Then, $(u^{j+1}_0 , v^{j+1}_0) = (u^j_0 + c_0^j, v_j^0 + d_0^j)$ in turn fulfill \eqref{eq:proof:CC-for-u(j)}.
		      Moreover, the estimates above, \eqref{eq:proof:CC-for-u(j)-3}, and the usual Sobolev trace theorem yield
		      \begin{align}\label{eq:est-corrections}
			      \notag & \smallnorm{(u_0,v_0) - (u_0^{j+1}, v_0^{j+1})}_{H_\Gamma^{m}(\Omega)} \leq \smallnorm{(c^j_0,d^j_0 )}_{H^m(\omega)} +\smallnorm{(u_0,v_0) - (u_0^j, v_0^j)}_{H_\Gamma^{m}(\Omega)}               \\
			      \notag & \lesssim \smallnorm{\tr \partial_\nu^j(c^j_0,d^j_0 )}_{H^{m-j-1/2}(\partial\omega)}+\smallnorm{(u_0,v_0) - (u_0^j, v_0^j)}_{H_\Gamma^{m}(\Omega)}                                                \\
			      \notag & \lesssim \norm{\big\llbracket( u^{j}_{(j)} \times \nu^\Gamma, v^{j}_{(j)} \times \nu^\Gamma)\big\rrbracket - (\pt^j \kappa_n(0), \pt^j \varkappa_n(0))}_{H^{m-j-1/2}(\Gamma)}                    \\
			      \notag & \seq + \smallnorm{\tr_t u^{j}_{(j)} - \pt^j \bndInh_n(0)}_{H^{m-j-1/2}(\partial\Omega)} +\smallnorm{(u_0,v_0) - (u_0^j, v_0^j)}_{H_\Gamma^{m}(\Omega)}                                           \\
			      \notag & \lesssim \norm{\big\llbracket (u^{\text{app}}_{(j)} \times \nu^\Gamma , v^{\text{app}}_{(j)} \times \nu^\Gamma)\big\rrbracket - (\pt^j \kappa_n(0), \pt^j \varkappa_n(0))}_{H^{m-j-1/2}(\Gamma)} \\
			      \notag & \seq+ \smallnorm{\tr_t u^{\text{app}}_{(j)} - \pt^j \bndInh_n(0)}_{H^{m-j-1/2}(\partial\Omega)} +\smallnorm{(u_0,v_0) - (u_0^j, v_0^j)}_{H_\Gamma^{m}(\Omega)}                                   \\
			             & \leq \hat C r,
		      \end{align}
		      where we also used the assumption and induction hypothesis.

		      Iteratively, this allows us to construct
		      $$
			      (u_0^m,v_0^m) = (u_0, v_0) + \sum_{j=1}^{m-1} (c_0^j, d_0^j) \in H^{m+a}_\Gamma(\Omega),
		      $$
		      satisfying the compatibility condition $(CC_m)$ of order $m$ and \eqref{eq:approx-close}.

		\item In order to pass to $(CC_{m+a})$, for $a>0$, we can construct $(\hat c^{m+a-1}_0,\hat d^{m+a-1}_0)$ correcting the compatibility condition iteratively as above.
		      However, we cannot control the norm through the (approximate) compatibility condition of $(u_0,v_0)$ as in \eqref{eq:est-corrections} before.
		      Instead, we rely on a scaling argument in this case.
		      In the following we let the $a=1$ in the compatibility condition, i.e., we show how to obtain $(CC_{m+1})$.
		      We will only focus on $\hat c^m_0$ as the argument for $\hat d^m_0$ follows similarly.

		      We take a covering $\{U_i\}_{i=0}^M$ of $\omega$ and coordinate charts $\{(U_i, \Psi_i)\}_{i=1}^M$ at the boundary of $\omega$ such that $(\cN_b^{\partial \Omega} \cup \cN^{\Gamma}_b )\cap \omega \subseteq \cup_{i=1}^M U_i$, $\Psi_i(U_i \cap \omega) = V_i \cap \mathbb{R}^3_+$ and $\Psi_i(U_i \cap \partial \omega) = V_i \cap \{y \in \R^3 \mid y_3 = 0\}$ for $i\geq 1$.
		      Here \[\cN_b^\Xi \coloneqq \{ x \in \overline\Omega \mid \mathrm{dist}(x, \Xi)<b\}\] denotes collars of width $b>0$ around $\Xi \in \{\partial \Omega, \Gamma\}$.
		      Let $\phi_i$ be a partition of unity subordinate to $\{U_i\mid i = 0, \dots, M\}$.

		      Let $\delta \in (0,1)$.
		      For $i=1, \dots, M$, let
		      \[
			      p_i \colon \R^3_+ \to \R^3; \qquad p_i(y) \coloneqq
			      \begin{cases}
				      (\phi_i \hat c_0^m)(\Psi_i^{-1}(y)), \quad & \text{for } y \in V_i     \\
				      0, \quad                                   & \text{for } y \notin V_i,
			      \end{cases}
		      \]
		      denote the local pushforward of $\phi_i \hat c_0^m$ extended by zero. We use the scaling matrix
		      $S_\delta \coloneqq \begin{psmallmatrix}
				      1&0&0\\
				      0&1&0\\
				      0&0&\tfrac{1}{\delta}
			      \end{psmallmatrix}$
		      and set
		      \[
			      c^{m,\delta}_0(x) \coloneqq \delta^m \sum_{i=1}^M p_i\left(S_\delta\Psi_i(x)\right), \qquad x \in \omega,
		      \]
		      where $p_i(S_\delta\Psi_i(x))$ is understood to be $0$ for $x \notin U_i$.
		      Thus, by substitution
		      \begin{align*}
			      \smallnorm{c^{m,\delta}_0}^2 _{L^2(\Omega)}
			       & \lesssim \sum_{i=1}^M \int_{\mathbb{R}^3_+} \left| \delta^m p_i\left(y', \frac{y_3}{\delta}\right) \right|^2 \dd y \\
			       & = \sum_{i=1}^M \delta^{2m} \int_{\mathbb{R}^3_+} \left| p_i(u) \right|^2 \delta \dd u                              \\
			       & \lesssim \delta^{2m+1} \smallnorm{\hat c_0^{m}}_{L^2(\Omega)}^2.
		      \end{align*}
		      Likewise, for the Sobolev norm of order $m$, we consider derivatives $D^\alpha$ with $|\alpha| \leq m$. The dominant terms in $\delta$ arise from the normal derivatives $\partial_{x_3}^\ell$ where $\ell = \alpha_3 \leq m$. Namely, for $\alpha = (\alpha', \ell)$ with $\ell \in \{1, \dots, m\}$, we get
		      \begin{align*}
			      \smallnorm{D^\alpha c^{m,\delta}_0}_{L^2(\Omega)}^2
			       & \lesssim \sum_{i=1}^M \int_{\mathbb{R}^3_+} \left| \delta^m \delta^{-\ell} (\partial_{y'}^{\alpha'} \partial_{y_3}^\ell p_i)\left(y', \frac{y_3}{\delta}\right) \right|^2 \dd y \\
			       & = \sum_{i=1}^M \delta^{2(m-\ell)} \int_{\mathbb{R}^3_+} \left| (D^\alpha p_i)(u) \right|^2 \delta \dd u                                                                         \\
			       & \leq C\delta^{2(m-\ell)+1} \smallnorm{\hat c_0^{m}}_{H_\Gamma^{\abs\alpha}(\Omega)}^2
		      \end{align*}
		      for some constant $C>0$.
		      Since $\ell \leq m$, the exponent $2(m-\ell)+1$ is larger or equal to $1$. Thus, for $\delta \to 0$ we have $\smallnorm{c^{m,\delta}_0}_{H^m(\omega)} \to 0$. Without loss of generality we may therefore assume that $\norm{c^{m}_0}_{H^m(\omega)} \leq r$ by choosing $c^{m}_0 \coloneqq c_0^{m, \delta}$ with $\delta \coloneqq \tfrac{r^2}{C\smallnorm{\hat c_0^{m}}_{H^m(\Omega)}^2}$.

		      Note that the $m$-th order normal derivative at $\partial\omega = \partial\Omega \cupdot \Gamma$ is unaffected by the scaling.
		      To see this, let us consider the partial derivatives of order $m$.
		      The derivatives of $(\phi_i \hat c_0^{m})\left(\Psi_i^{-1}\left(S_\delta\Psi_i(x)\right)\right)$ involving derivatives of $\hat c_0^{m}$ of order strictly less than $m$ vanish on $\partial\omega $, since $\hat c_0^{m} \in H^m_0(\omega)$. Hence, for $\abs\alpha = m$
		      \[
			      D^\alpha c^{m,\delta}_0(x) = \delta^m \sum_{i=1}^M \phi_i(x) D^\alpha \left(\hat c_0^{m}\left(\Psi_i^{-1}\left(S_\delta\Psi_i(x)\right)\right)\right)
		      \]
		      on $\partial\omega$. In the following we are thus only interested in terms where `all derivatives fall upon $\hat c_0^{m}$'.

		      Let $p_{\delta,i}$ be the scaled function in local coordinates, i.e.,
		      \[
			      p_{\delta,i}\colon \R^3 \to \R^3; \qquad p_{\delta,i}(y)=\delta^m p_i(S_\delta y).
		      \]
		      Note that by the chain rule applied to the $m$-th derivative, only the normal derivatives (in local coordinates) scale with $\delta^{-1}$ by the definition of $S_\delta$. Specifically, on the boundary $\{y_3=0\}$
		      \[
			      \partial_{y_3}^m p_{\delta,i}(y) = \delta^m \left(\frac{1}{\delta}\right)^m (\partial_{y_3}^m p_i)(y) = (\partial_{y_3}^m p_i)(y).
		      \]
		      By the chain rule we get
		      \[
			      D c^{m,\delta}_0(x) = \sum_{i=1}^M D|_{(\Psi_i(x))} p_{\delta,i} \cdot D \Psi_i(x) = \sum_{i=1}^M(\partial_{y_3} p_{\delta,i})|_{y=\Psi_i(x)} \left(\nabla (\Psi_i)_3\right)^\top \quad \text{on }\partial\omega,
		      \]
		      since $p_{\delta,i}|_{\{y_3=0\}} = 0$ and hence $D p_{\delta,i} = (0 \mid 0 \mid \partial_{y_3} p_{\delta,i})$ on $\{y_3=0\}$. Here $(\cdot)_3$ denotes the third component.
		      Thus, for $x \in \partial\omega$
		      \begin{equation*}
			      \partial_\nu c^{m,\delta}_0(x) = (\nu \cdot \nabla) c^{m,\delta}_0(x) = \sum_{i=1}^M (\nu\cdot\nabla (\Psi_i)_3) \big(\partial_{y_3} p_{\delta,i}\big)|_{y=\Psi_i(x)},
		      \end{equation*}
		      i.e., the normal derivative corresponds to $\partial_{y_3}$ in local coordinates times the factor $(\nu\cdot\nabla (\Psi_i)_3)$ depending on the chart.
		      (Since $\partial\omega\cap U_i =\{x \in U_i\mid (\Psi_i(x))_3=0\}$ the gradient $\nabla (\Psi_i)_3$ is normal to $\partial\omega$.)
		      Similarly, for the $m$-th derivative we obtain
		      \begin{align*}
			      \partial_\nu^m c^{m,\delta}_0(x) & = \sum_{i=1}^M (\nu\cdot\nabla (\Psi_i)_3)^m \big(\partial_{y_3}^m p_{\delta,i}\big)|_{y=\Psi_i(x)}                                                     \\
			                                       & = \sum_{i=1}^M \phi_i(x) (\nu\cdot\nabla (\Psi_i)_3)^m \partial_{y_3}^m\left(\hat c_0^{m}\left(\Psi_i^{-1}\left( y \right)\right)\right)|_{y=\Psi_i(x)} \\
			                                       & = \partial_\nu^m \hat c_0^{m}(x).
		      \end{align*}
		      Note that here we can again neglect terms where derivatives fall upon $\phi_i(x) (\nu\cdot\nabla (\Psi_i)_3)$ and also that $D^\beta p_{\delta,i}|_{\{y_3=0\}} = 0$ for $\abs{\beta} \leq m-1$ because $\hat c_0^{m} \in H^m_0(\omega)$.

		      For $d^m_0$ we repeat the same argument as for $c^m_0$.
		      Thus, $(u_0^{m+1}, v_0^{m+1}) \coloneqq (u_0^{m}, v_0^{m}) + (c^m_0,d^m_0 )$ satisfies $(CC_{\alpha_n, \beta_n, \sigma_n, \rho_n, f_n, g_n, \bndInh_n, \kappa_n, \varkappa_n, m+1})$. Iterating this procedure for $\hat \imath \leq a$ yields $(u_0^{m+a}, v_0^{m+a})$ satisfying $(CC_{\alpha_n, \beta_n, \sigma_n, \rho_n, f_n, g_n, \bndInh_n, \kappa_n, \varkappa_n, m+a})$.
		      \qedhere
	\end{enumerate}
\end{proof}

We used the following lemma in the proof above.

\begin{lemma}[Simplification of $(\alpha \curl)^m$ for vanishing traces]\label{lem:Appendix:tr-curl}
	Let $m \in \N$ and let $\alpha_i \in C^\infty(\overline{\omega}; \Rsym)$ be positive definite matrices.
	Let $u \in (H^{m+1}(\omega) \cap H^m_0(\omega))^3$.
	Define the operator $L_i u \coloneqq \alpha_i (\curl u)$.
	Then
	\[
		\tr(\Pi_{i=1}^m L_i u) = \left(\Pi^m_{i=1} T_{\alpha_i, x}\right) \tr(\partial_\nu^m u) \qquad \text{on } \partial\omega,
	\]
	where $T_{\alpha_i, x}\colon \R^3 \to \R^3$ is the linear map defined for $x \in \partial \omega$ by
	\begin{equation}
		\label{eq:def-T_alpha,x}
		T_{\alpha_i, x}(\xi) \coloneqq \alpha_i(x) (\nu(x) \times \xi).
	\end{equation}
\end{lemma}

\begin{proof}
	\begin{enumerate}[label = \arabic*), wide, labelwidth=0pt, labelindent=0pt, itemsep=0.5em, topsep=0pt]
		\item \label{step1-curl} We first consider the case $m=1$.
		      Let $u \in (H^2(\omega) \cap H^1_0(\omega))^3$ and let $u_n \in C^\infty(\overline \omega)^3$ with $\tr u_n =0$ approximating $u$ in $H^2(\omega)$.
		      Since $\tr u_n =0$, the gradient runs orthogonal to $\partial \omega$, i.e.,
		      \[
			      \partial_j (u_n)_k = (\nabla (u_n)_k)_j = (\nu \cdot \nabla (u_n)_k) \nu_j \qquad \text{for }j,k \in\{1,2,3\}.
		      \]
		      By approximation this equality also holds for $u$.
		      The curl now becomes
		      \[
			      (\curl u)_i = \varepsilon_{ijk} \partial_j u_k = \varepsilon_{ijk}\nu_j(\nu \cdot \nabla u_k) = (\nu \times \partial_\nu u)_i,
		      \]
		      where $\varepsilon_{ijk}$ for $i,j,k \in\{1,2,3\}$ denotes the Levi-Civita symbol. Thus,
		      \[
			      \tr(\alpha_1 \curl u) = T_{\alpha_1, x} \partial_\nu u.
		      \]

		\item Now let us consider the general case $u \in H^{m+1} \cap H^m_0$.
		      We aim to compute the trace of $\Pi_{i=1}^m L_i u = \Pi_{i=1}^m(\alpha_i \curl) u$.
		      For $u \in H^m_0(\omega)$ the normal derivatives up to order $m-1$ are zero
		      \[
			      \partial_\nu^\ell u |_{\partial \omega} = 0 \quad \text{for } \ell = 0, 1, \dots, m-1.
		      \]
		      By the product rule
		      \[
			      \Pi_{i=1}^m L_i u = P^m u + R(u, Du, \dots, D^{m-1}u),
		      \]
		      where $P^m$ is the principal part of $\Pi_{i=1}^m L_i$ and $R$ is a linear lower order differential operator. Since $u \in H^{m}_0(\omega)$ the remainder $R(u, Du, \dots, D^{m-1}u)$ vanishes on $\partial \omega$.

		      So, an iteration of Step~\ref{step1-curl} yields
		      \[
			      \tr(\Pi_{i=1}^m L_i u) = T_{\alpha_m, x} (\partial_\nu ( T_{\alpha_{m-1}, x} (\partial_\nu (\dots T_{\alpha_1, x} (\partial_\nu u)\dots)))).
		      \]
		      Moreover, the commutator $[\partial_\nu,T_{\alpha_i, x}]$ results in lower order terms, which again vanish on the boundary, i.e.,
		      \[
			      \partial_\nu (T_{\alpha_i, x} \psi) = (\partial_\nu T_{\alpha_i, x}) \psi + T_{\alpha_i, x} (\partial_\nu \psi) = T_{\alpha_i, x} (\partial_\nu \psi).
		      \]
		      Thus, we commute all normal derivatives to the right obtaining
		      \[
			      \tr(\Pi_{i=1}^m L_i u) = (T_{\alpha_m,x} \circ \dots \circ T_{\alpha_1,x}) \partial_\nu^m u \qquad \text{on } \partial \omega. \qedhere
		      \]
	\end{enumerate}
\end{proof}

\section{Convergence Result}
\label{sec:3-convergence-result}

Based on the previous section, we will work in the following setting:
\smallskip
\begin{enumerate}[leftmargin=25pt, label=\textup{($\text{H}_\text{a}$)}]
	\item \label{hypothesis-approx}\quad \parbox[m]{\linegoal}{
		      Let $a \in \N_0$, $m \in \N$ with $\tilde m + 2 +a\leq k$.
		      Let $(u,v) \in \cG_\Gamma^m$ be a solution to \eqref{eq:maxwell-lem}--\eqref{eq:maxwell-lem4} with initial values $(u_0,v_0)$ and coefficients satisfying \eqref{eq:maxwell-lem-coeff} -- \eqref{eq:maxwell-lem-IVs}.
		      We approximate our coefficients as in \ref{def:approx-data}. By Lemma~\ref{lem:approx-CC1} and Proposition~\ref{prop:fix-CC} we may choose initial values $(u^{m+a}_{0,n}, v^{m+a}_{0,n})\in H^{m+a}_\Gamma(\Omega)$ satisfying $(CC_{\alpha_n, \beta_n, \sigma_n, \rho_n, f_n, g_n, \bndInh_n, \kappa_n, \varkappa_n, m+a})$ converging to the original initial data $(u_0,v_0) \in H^m_\Gamma(\Omega)$. We denote the corresponding solutions by $(u_n,v_n) \in \cG_\Gamma^{m+a}$.}
\end{enumerate}

\medskip

We want to show the convergence to our original solution $(u,v)$ in $\cG_\Gamma^m$. Unfortunately, considering the difference between the approximate solution and the exact solution
\begin{equation*}
	U \coloneqq u_{n} - u \quad \text{and} \quad V \coloneqq v_{n} - v
\end{equation*}
leads to the system
\begin{alignat*}{-1}
	\begin{split}
		\alpha \pt U & = \hphantom{-}\curl V - \sigma U \hspace{-0.5pt}- \mathfrak{f}_n, \\
		\beta \pt V  & = - \curl U - \rho V \hspace{0.5pt}- \mathfrak{g}_n,
	\end{split} \qquad & t                                                                                                     & \geq 0, x \in \Omega_\pm,                                                                                                   \\
	\notag U \times \nu                                                                 & = \bndInh_n - \bndInh, \qquad                                                                         & t                         & \geq 0, x \in \partial \Omega, \\
	\notag \llbracket U \times \nu^\Gamma \rrbracket                                    & = \kappa_n - \kappa, \quad \llbracket V \times \nu^\Gamma \rrbracket = \varkappa_n - \varkappa,\qquad & t                         & \geq 0, x \in \Gamma,          \\
	\notag U(0)                                                                         & = u_{0,n} - u_0, \quad V(0) = v_{0,n} - v_0,\qquad                                                    & x                         & \in \Omega_\pm.
\end{alignat*}
The inhomogeneities $\mathfrak{f}_n$ and $\mathfrak{g}_n$ lie in $L^2(0,T; L^2(\Omega)^3)$ and are given by
\begin{alignat*}{-1}
	\mathfrak{f}_{n} & \coloneqq (f_n - f) &  & - (\alpha - \alpha_n)\pt u_{n} &  & + (\sigma - \sigma_n) u_{n}, \\
	\mathfrak{g}_{n} & \coloneqq (g_n - g) &  & - (\beta - \beta_n)\pt v_{n}   &  & + (\rho - \rho_n) v_{n}.
\end{alignat*}
In this direct approach the terms $(\alpha - \alpha_n)\pt u_{n}$ and $(\beta - \beta_n)\pt v_{n}$ cause trouble. While the difference of the coefficients vanishes in the right norm, one lacks control over the time derivative.
However, the a-priori estimate \eqref{eq:lwp-apriori} at least yields convergence in $\cG_\Gamma^{m-1}$.

To obtain $\cG_\Gamma^m$-convergence, we first mollify the original solution in time direction.
This allows us to show that the time derivatives converge. Afterwards, we inductively reduce the convergence of spatial norms to the temporal ones using a div-curl lemma, see Theorem~\ref{thm:div-curl-perfect}, and the Maxwell system.

\begin{definition}[Mollifier]\label{def:mollifier}
	In the following, let $\eta_\epsilon(t) = \frac{1}{\epsilon}\eta\left(\frac{t}{\epsilon}\right)$ for the usual mollifier
	\[
		\eta(t) \coloneqq \begin{cases}
			C \exp\left(\frac{1}{\abs{t}^2-1}\right), \quad & \abs{t} \leq 1, \\
			0, \quad                                        & \abs{t} >1,
		\end{cases}
	\]
	where $C>0$ is chosen such that $\int_{\mathbb{R}} \eta(t) \dd t = 1$. Let $(u)_\epsilon \coloneqq u * \eta_\epsilon$ denote the mollification in time.
\end{definition}

We will need the following elementary mollification result later on.

\begin{lemma}\label{lem:commutator-mollifier}
	Let $T > 0$, $\epsilon_0 > 0$, $u \in C\bigl([-\epsilon_0, T+\epsilon_0]; L^2(\Omega)^3\bigr)$ and $\alpha \in C^1\bigl([-\epsilon_0, T+\epsilon_0]; H^2_\Gamma(\Omega, \Rsym)\bigr)$. Then,
	\[
		\pt (\alpha u)_\epsilon - \alpha \pt u_\epsilon \to (\pt \alpha) u \quad \text{ in } C\bigl([0, T]; L^2(\Omega)^3\bigr),
	\]
	for $\epsilon< \epsilon_0$ and $\epsilon \to 0$.
\end{lemma}

\begin{proof}
	In the following, let $\epsilon < \epsilon_0$.
	Let $\cM_\epsilon(t) = \pt (\alpha u)_\epsilon(t) - \alpha(t) \pt u_\epsilon(t)$. Expanding the convolutions and applying the change of variables $\tau = \frac{t-s}{\epsilon}$, the commutator evaluated at time $t \in [0,T]$ is given by
	\[
		\cM_\epsilon(t) = \int_{-1}^{1} \frac{\alpha(t-\epsilon\tau) - \alpha(t)}{\epsilon} u(t-\epsilon\tau) \eta'(\tau)\dd \tau.
	\]
	We rewrite the limit as
	\[
		\pt \alpha(t) u(t) = -\pt \alpha(t) u(t) \int_{-1}^{1} \tau \eta'(\tau)\dd \tau
	\]
	and consider the difference
	\[
		\cM_\epsilon(t) - (\pt \alpha(t)) u(t) = I_1(t, \epsilon) + I_2(t, \epsilon),
	\]
	where
	\begin{align*}
		I_1(t, \epsilon) & = \int_{-1}^{1} \left[ \frac{\alpha(t-\epsilon\tau) - \alpha(t)}{\epsilon} + \tau \pt \alpha(t) \right] u(t-\epsilon\tau) \eta'(\tau)\dd \tau, \\
		I_2(t, \epsilon) & = - \int_{-1}^{1} \tau \pt \alpha(t) \big( u(t-\epsilon\tau) - u(t) \big) \eta'(\tau)\dd \tau.
	\end{align*}

	Since 
	$\alpha \in C^1([-\epsilon_0, T+\epsilon_0]; L^\infty(\Omega)^{3\times 3})$ by the Sobolev embeddings, we see that
	\begin{align*}
		\sup_{t \in [0,T], \tau \in [-1,1]} \left\| \frac{\alpha(t-\epsilon\tau) - \alpha(t)}{\epsilon} + \tau \pt \alpha(t) \right\|_{L^\infty(\Omega)} & \to 0, \\
		\sup_{t \in [0,T], \tau \in [-1,1]} \|u(t-\epsilon\tau) - u(t)\|_{L^2(\Omega)}                                                                   & \to 0,
	\end{align*}
	for $\epsilon \to 0$. Hence,
	\[
		\sup_{t \in [0,T]} I_1(t, \epsilon), \sup_{t \in [0,T]} I_2(t, \epsilon) \to 0 \qquad \text{for }\epsilon \to 0
	\]
	and the claim follows.
\end{proof}

Recall that the $\cG_\Gamma^m$ norm is given by
\[
	\norm{u}_{\cG_\Gamma^m([0,T] \times \Omega)} \coloneqq \sup_{t \in [0,T]} \sum_{j=0}^m \smallnorm{\pt^j u(t)}_{H^{m-j}_\Gamma}.
\]
We first focus on the highest time derivative, i.e., $\sup_{t \in [0,T]} \smallnorm{\pt^m u(t)}_{L^2(\Omega)}$.

\begin{lemma}[Convergence of Time Derivatives]\label{lem:time-convergence}
	Let \ref{hypothesis-approx} hold. Then,
	\[
		\lim_{n \to \infty} \sup_{t \in [0,T]} \Big( \norm{\pt^m u_n(t) - \pt^m u(t)}_{L^2(\Omega)} + \norm{\pt^m v_n(t) - \pt^m v(t)}_{L^2(\Omega)} \Big) = 0.
	\]
\end{lemma}

\begin{proof}
	\begin{enumerate}[label = \arabic*), wide, labelwidth=0pt, labelindent=0pt, itemsep=0.5em, topsep=0pt]
		\item In order to apply a mollifier on $[0,T]$, we continuously extend the coefficients $\alpha, \alpha_n, \beta, \beta_n \in \cC_\eta^{\tilde m}$, $\sigma, \sigma_n, \rho, \rho_n \in \cC^{\tilde m}$, $f, f_n, g, g_n \in H^m([0,T] \times \Omega_\pm)$, and the extensions of the boundary and interface data $\bndInh^\circ,\bndInh_n^\circ ,\kappa^\circ,\kappa_n^\circ,\varkappa^\circ,\varkappa_n^\circ \in \cG^{m+1}_\Gamma$ to a slightly larger time interval $[-\epsilon_0, T+ \epsilon_0]$ by the usual extension operator as in Theorem~5.19 in \cite{adams2003}.
		      Theorem~3.1 in \cite{schnaubelt2022} then provides extensions of the solutions $(u,v), (u_n,v_n)$ in $\cG_\Gamma^{m}([0,T])$ respectively $\cG_\Gamma^{m + a}([0,T])$ to $\cG_\Gamma^{m}([-\epsilon_0, T+ \epsilon_0])$ respectively $\cG_\Gamma^{m + a}([-\epsilon_0, T+ \epsilon_0])$.
		      From now on, we denote the extensions by the same letters.

		      Let $U_n \coloneqq \pt^m u_n$, $V_n \coloneqq \pt^m v_n$ and $U \coloneqq \pt^m u$, $V \coloneqq \pt^m v$ denote the highest-order time derivatives.
		      We mollify $U_n$ on $[0,T]$, yielding $U_{n,\delta} \coloneqq (U_n)_\delta$ for $0<\delta<\epsilon_0$. Likewise, we mollify $U$ on $[0,T]$ and obtain $U_\epsilon \coloneqq (U)_\epsilon$ for $0<\epsilon<\epsilon_0$.
		      (For the approximating sequence this has to be done to include the case of $a=0$. Otherwise, the proof simplifies and only requires the mollification of $(U,V)$. Regardless, the additional arguments are quite similar to that case.)

		      We define the differences $D^u_{n,\delta,\epsilon} \coloneqq U_{n,\delta} - U_\epsilon$ and $D^v_{n,\delta,\epsilon} \coloneqq V_{n,\delta} - V_\epsilon$, and set
		      \begin{equation*}
			      e^D_{n,\delta,\epsilon}(t) \coloneqq \frac{1}{2} \int_\Omega \Big( D^u_{n,\delta,\epsilon} \cdot \alpha_n D^u_{n,\delta,\epsilon} + D^v_{n,\delta,\epsilon} \cdot \beta_n D^v_{n,\delta,\epsilon} \Big) \dd x.
		      \end{equation*}

		      Moreover, as in \eqref{eq:simultaneous-extensions} there exist extensions $(W^u_n, W^v_n)$ and $(W^u, W^v)$ in $\cG^{m+1}_\Gamma$ matching the data simultaneously, i.e.,
		      \begin{alignat*}{-1}
			      W^u \times \nu   & = \bndInh   &  & \text{ on }\partial\Omega, & \qquad \llbracket W^u \times \nu^\Gamma \rrbracket   & = \kappa,   & \qquad \llbracket W^v \times \nu^\Gamma \rrbracket                   & = \varkappa,   \\*
			      W^u_n \times \nu & = \bndInh_n &  & \text{ on }\partial\Omega, & \qquad \llbracket W^u_n \times \nu^\Gamma \rrbracket & = \kappa_n, & \quad \text{and} \quad \llbracket W^v_n \times \nu^\Gamma \rrbracket & = \varkappa_n.
		      \end{alignat*}
		      Because the extensions lie in $\cG^{m+1}_\Gamma$, their $m$-th time derivatives inherently belong to $H^1_\Gamma(\Omega) \subset H(\curl,\Omega_\pm)$ with the $(m+1)$-th time derivatives in $L^2(\Omega)$.
		      Furthermore,
		      \begin{align*}
			      \sup_{t \in [-\epsilon_0,T+\epsilon_0]} & \left(\norm{\pt^mW^u_n - \pt^mW^u}_{H^1_\Gamma(\Omega)} + \norm{\pt^mW^v_n - \pt^mW^v}_{H^1_\Gamma(\Omega)} \right)                    \\
			                                              & \lesssim \norm{(\bndInh_n^\circ-\bndInh^\circ, \kappa_n^\circ - \kappa ,\varkappa_n^\circ - \varkappa^\circ)}_{\cG^{m+1}_\Gamma} \to 0
		      \end{align*}
		      for $n \to \infty$.

		      Let $E^u_{n,\delta,\epsilon} \coloneqq (\pt^m W^u_n)_\delta - (\pt^m W^u)_\epsilon$ and $E^v_{n,\delta,\epsilon} \coloneqq (\pt^m W^v_n)_\delta - (\pt^m W^v)_\epsilon$.
		      We define the shifted differences
		      \[
			      \tilde{D}^u_{n,\delta,\epsilon} \coloneqq D^u_{n,\delta,\epsilon} - E^u_{n,\delta,\epsilon} \quad \text{and} \quad \tilde{D}^v_{n,\delta,\epsilon} \coloneqq D^v_{n,\delta,\epsilon} - E^v_{n,\delta,\epsilon},
		      \]
		      which by construction satisfy homogeneous boundary and jump conditions, allowing us to integrate by parts.

		      We set the shifted energy
		      \begin{equation*}
			      \tilde{e}^D_{n,\delta,\epsilon}(t) \coloneqq \frac{1}{2} \int_\Omega \Big( \tilde{D}^u_{n,\delta,\epsilon} \cdot \alpha_n \tilde{D}^u_{n,\delta,\epsilon} + \tilde{D}^v_{n,\delta,\epsilon} \cdot \beta_n \tilde{D}^v_{n,\delta,\epsilon} \Big) \dd x.
		      \end{equation*}
		      Applying the time derivative in the classical sense yields
		      \begin{equation}
			      \label{eq:energy_deriv}
			      \begin{aligned}
				      \frac{\dd}{\dd t} \tilde{e}^D_{n,\delta,\epsilon}(t) & = \int_\Omega \Big( \tilde{D}^u_{n,\delta,\epsilon} \cdot \alpha_n \pt \tilde{D}^u_{n,\delta,\epsilon} + \tilde{D}^v_{n,\delta,\epsilon} \cdot \beta_n \pt \tilde{D}^v_{n,\delta,\epsilon} \Big) \dd x                     \\
				                                                           & \seq + \frac{1}{2} \int_\Omega \Big( \tilde{D}^u_{n,\delta,\epsilon} \cdot (\pt \alpha_n) \tilde{D}^u_{n,\delta,\epsilon} +\tilde{D}^v_{n,\delta,\epsilon} \cdot (\pt \beta_n) \tilde{D}^v_{n,\delta,\epsilon}\Big) \dd x.
			      \end{aligned}
		      \end{equation}

		      In the following, we will only consider $u$ and $u_n$, since the calculations for $v$ and $v_n$ are done similarly.
		      We intend to apply a Gronwall argument to \eqref{eq:energy_deriv}. For this, we first expand $\alpha_n \pt \tilde{D}^u_{n,\delta,\epsilon} = \alpha_n \pt D^u_{n,\delta,\epsilon} - \alpha_n \pt E^u_{n,\delta,\epsilon}$ to substitute expressions for the time derivatives.

		      The $(m-1)$-th level system for $u_n$ is
		      \begin{equation}
			      \label{eq:not-mollifed}
			      \alpha_n U_n = \curl(\pt^{m-1} v_n) - \sigma_n \pt^{m-1} u_n - \pt^{m-1} f_n - \cD^u_n,
		      \end{equation}
		      where $\cD^u_n \coloneqq [\pt^{m-1}, \alpha_n \pt] u_n + [\pt^{m-1}, \sigma_n] u_n$ is the commutator with the derivative.
		      We further set $\cM^n_\delta \coloneqq \pt(\alpha_n U_n)_\delta - \alpha_n \pt U_{n,\delta}$, mollify \eqref{eq:not-mollifed} and take the time derivative.
		      It follows that
		      \begin{equation}\label{eq:ME-for-pt^m-u_n}
			      \alpha_n \pt U_{n,\delta} = \curl V_{n,\delta} - \pt(\sigma_n \pt^{m-1}u_n)_\delta - (\pt^m f_n)_\delta - \pt(\cD^u_n)_\delta - \cM^n_\delta.
		      \end{equation}

		      Similarly, for $U$ we mollify the $(m-1)$-th level equation
		      \[
			      \alpha U = \curl (\pt^{m-1} v) - \sigma \pt^{m-1} u - \pt^{m-1} f - \cD^u
		      \]
		      with $\cD^u \coloneqq [\pt^{m-1}, \alpha \pt] u + [\pt^{m-1}, \sigma] u$. Also applying the time derivative and introducing $\cM_\epsilon \coloneqq \pt (\alpha U)_\epsilon - \alpha \pt U_\epsilon$ yields
		      \begin{equation}\label{eq:ME-for-pt^m-u}
			      \alpha \pt U_\epsilon = \curl V_\epsilon - \pt(\sigma \pt^{m-1} u)_\epsilon - (\pt^m f)_\epsilon - \pt (\cD^u)_\epsilon - \cM_\epsilon.
		      \end{equation}
		      Writing $\alpha_n \pt U_\epsilon = \alpha \pt U_\epsilon + (\alpha_n - \alpha)\pt U_\epsilon$, Formulas~\eqref{eq:ME-for-pt^m-u_n} and \eqref{eq:ME-for-pt^m-u} provide the equation
		      \[
			      \alpha_n \pt D^u_{n,\delta,\epsilon} = \curl D^v_{n,\delta,\epsilon} - (\sigma_n+m \pt \alpha_n) D^u_{n,\delta,\epsilon} + R^u_{n,\delta,\epsilon}
		      \]
		      for the unshifted difference, where the remainder $R^u_{n,\delta,\epsilon}$ is given by
		      \begin{align}\label{eq:proof-remainder-terms-unshifted}
			      \notag R^u_{n,\delta,\epsilon} & \coloneqq - (\alpha_n - \alpha)\pt U_\epsilon+ \sigma_n U_{n,\delta} - \sigma_n U_\epsilon+ m(\pt \alpha_n) U_{n,\delta} - m(\pt \alpha_n) U_\epsilon \\
			                                     & \seq + \pt(\sigma \pt^{m-1} u)_\epsilon - \pt(\sigma_n \pt^{m-1} u_n)_\delta + (\pt^m f)_\epsilon - (\pt^m f_n)_\delta                                \\
			      \notag                         & \seq + \pt(\cD^u)_\epsilon - \pt(\cD^u_n)_\delta + \cM_\epsilon - \cM^n_\delta .
		      \end{align}
		      We added the $- (\sigma_n+m \pt \alpha_n) D^u_{n,\delta,\epsilon}$-term to rebalance the remainder later on.

		      Substituting the shifted variables, $(\tilde{D}^u_{n,\delta,\epsilon}, \tilde{D}^v_{n,\delta,\epsilon})$ satisfy
		      \[
			      \alpha_n \pt \tilde{D}^u_{n,\delta,\epsilon} = \curl \tilde{D}^v_{n,\delta,\epsilon} - (\sigma_n+m \pt \alpha_n) \tilde{D}^u_{n,\delta,\epsilon} + \tilde{R}^u_{n,\delta,\epsilon},
		      \]
		      with the modified remainder
		      \begin{equation}
			      \label{eq:proof-remainder-terms-shifted}
			      \tilde{R}^u_{n,\delta,\epsilon} \coloneqq R^u_{n,\delta,\epsilon}- \alpha_n \pt E^u_{n,\delta,\epsilon} + \curl E^v_{n,\delta,\epsilon} - (\sigma_n + m\pt \alpha_n) E^u_{n,\delta,\epsilon}.
		      \end{equation}
		      Analogously, the difference of the magnetic terms satisfies
		      \[
			      \beta_n \pt \tilde{D}^v_{n,\delta,\epsilon} = -\curl \tilde{D}^u_{n,\delta,\epsilon} - (\rho_n+ m \pt \beta_n)\tilde{D}^v_{n,\delta,\epsilon} + \tilde{R}^v_{n,\delta,\epsilon},
		      \]
		      where $\tilde{R}^v_{n,\delta,\epsilon} \coloneqq R^v_{n,\delta,\epsilon} - \curl E^u_{n,\delta,\epsilon} - \beta_n \pt E^v_{n,\delta,\epsilon} - (\rho_n + m\pt \beta_n) E^v_{n,\delta,\epsilon}$.

		      We substitute these equations into \eqref{eq:energy_deriv} and integrate by parts. Since $\tilde{D}^u_{n,\delta,\epsilon}$ and $\tilde{D}^v_{n,\delta,\epsilon}$ satisfy homogeneous boundary and jump conditions
		      \[
			      \tilde{D}^u_{n,\delta,\epsilon}\times \nu = 0 \text{ on }\partial\Omega \qquad\text{and}\qquad \bigl\llbracket \bigl(\tilde{D}^u_{n,\delta,\epsilon} \times \nu^\Gamma, \tilde{D}^v_{n,\delta,\epsilon} \times \nu^\Gamma \bigr)\bigr\rrbracket = 0,
		      \]
		      all surface integrals vanish. It follows that
		      \begin{align*}
			      \frac{\dd}{\dd t} \tilde{e}^D_{n, \delta,\epsilon}(t) & = \int_\Omega \Big( \tilde{D}^u_{n, \delta,\epsilon} \cdot (\curl \tilde{D}^v_{n, \delta,\epsilon} - (\sigma_n + m\pt \alpha_n) \tilde{D}^u_{n, \delta,\epsilon} + \tilde{R}^u_{n, \delta,\epsilon})\Big) \dd x                                                                \\
			                                                            & \seq + \int_\Omega \Big(\tilde{D}^v_{n, \delta,\epsilon} \cdot (-\curl \tilde{D}^u_{n, \delta,\epsilon} - (\rho_n + m\pt \beta_n) \tilde{D}^v_{n, \delta,\epsilon} + \tilde{R}^v_{n, \delta,\epsilon})\Big) \dd x                                                              \\
			                                                            & \seq + \frac{1}{2} \int_\Omega \Big( \tilde{D}^u_{n, \delta,\epsilon} \cdot (\pt \alpha_n) \tilde{D}^u_{n, \delta,\epsilon} + \tilde{D}^v_{n, \delta,\epsilon} \cdot (\pt \beta_n) \tilde{D}^v_{n, \delta,\epsilon}\Big) \dd x                                                 \\
			                                                            & = \int_\Omega \Big( \tilde{D}^u_{n, \delta,\epsilon} \cdot (- \sigma_n \tilde{D}^u_{n, \delta,\epsilon} + \tilde{R}^u_{n, \delta,\epsilon}) + \tilde{D}^v_{n, \delta,\epsilon} \cdot (- \rho_n \tilde{D}^v_{n, \delta,\epsilon} + \tilde{R}^v_{n, \delta,\epsilon})\Big) \dd x \\
			                                                            & \seq + \!\int_\Omega \Big( \tilde{D}^u_{n, \delta,\epsilon} \cdot (\tfrac{1}{2}-m)(\pt \alpha_n) \tilde{D}^u_{n, \delta,\epsilon} + \tilde{D}^v_{n, \delta,\epsilon} \cdot (\tfrac{1}{2}-m)(\pt \beta_n) \tilde{D}^v_{n, \delta,\epsilon}\Big) \dd x                           \\
			                                                            & \leq C_\eta \left(1 + \norm{(\pt \alpha_n, \pt \beta_n,\sigma_n, \rho_n)(t)}_{L^\infty}\right) \tilde{e}^D_{n, \delta,\epsilon}(t)                                                                                                                                             \\
			                                                            & \seq + \smallnorm{\tilde{R}^u_{n, \delta,\epsilon}(t)}_{L^2(\Omega)}^2 + \smallnorm{\tilde{R}^v_{n, \delta,\epsilon}(t)}_{L^2(\Omega)}^2
		      \end{align*}
		      for some constant $C_\eta>0$.

		      Gronwall's inequality then yields
		      \begin{equation}
			      \label{eq:proof:gronwall-approx}
			      \tilde{e}^D_{n, \delta,\epsilon}(t) \leq \e^{C_G t} \left( \tilde{e}^D_{n, \delta,\epsilon}(0) + \int_0^t \Big( \smallnorm{\tilde{R}^u_{n, \delta,\epsilon}(s)}_{L^2(\Omega)}^2 + \smallnorm{\tilde{R}^v_{n, \delta,\epsilon}(s)}_{L^2(\Omega)}^2 \Big) \dd s \right)
		      \end{equation}
		      with the constant $C_G \coloneqq \sup_{n \in \N, t \in [0,T]} C_\eta \left(1 + \norm{(\pt \alpha_n, \pt \beta_n,\sigma_n, \rho_n)(t)}_{L^\infty}\right)>0$.

		\item \label{step:proof-time-convergence-delta-limit}
		      For fixed $n$ and $\epsilon$, we take the limit $\delta \to 0$.
		      We define the limits $E^u_{n,\epsilon} \coloneqq \pt^m W^u_n - (\pt^m W^u)_\epsilon$ and $E^v_{n,\epsilon} \coloneqq \pt^m W^v_n - (\pt^m W^v)_\epsilon$, as well as the shifted differences $\tilde{D}^u_{n,\epsilon} \coloneqq D^u_{n,\epsilon} - E^u_{n,\epsilon}$ and $\tilde{D}^v_{n,\epsilon} \coloneqq D^v_{n,\epsilon} - E^v_{n,\epsilon}$, alongside the corresponding energy $\tilde{e}^D_{n,\epsilon}(t) \coloneqq \frac{1}{2} \int_\Omega \big( \tilde{D}^u_{n,\epsilon} \cdot \alpha_n \tilde{D}^u_{n,\epsilon} + \tilde{D}^v_{n,\epsilon} \cdot \beta_n \tilde{D}^v_{n,\epsilon} \big) \dd x$.

		      Since $(u_n, v_n) \in \cG_\Gamma^m$, and hence $\pt^m u_n, \pt^m v_n \in C([-\epsilon_0,T+\epsilon_0]; L^2(\Omega))$, standard results on mollifiers yield the uniform convergence $U_{n,\delta} \to U_n$ and $V_{n,\delta} \to V_n$ in $C([0,T]; L^2(\Omega))$ for $\delta \to 0$.
		      Furthermore, $(\pt^m f_n)_\delta$ tends to $\pt^m f_n$ in $L^2([0,T]\times \Omega)$, $\pt(\sigma_n \pt^{m-1} u_n)_\delta$ to $\pt(\sigma_n \pt^{m-1} u_n)$ in $C([0,T],L^2(\Omega))$, $\pt(\cD^u_n)_\delta$ to $\pt\cD^u_n$ in $C([0,T],L^2(\Omega))$, and by Lemma~\ref{lem:commutator-mollifier}, the commutator $\cM^n_\delta$ to $(\pt \alpha_n)U_n$ in $C([0,T],L^2(\Omega))$ for $\delta \to 0$.

		      Additionally, since $\pt^m(W^u_n, W^v_n)$ and $\pt^m(W^u, W^v)$ belong to $\cG^{1}_\Gamma$ the mollified functions $E^u_{n,\delta,\epsilon} \to E^u_{n,\epsilon}$ and $E^v_{n,\delta,\epsilon} \to E^v_{n,\epsilon}$ converge in $\cG^{1}_\Gamma$ for $\delta \to 0$. In particular, this guarantees the pointwise convergence of the energy $\tilde{e}^D_{n,\delta,\epsilon}(t) \to \tilde{e}^D_{n,\epsilon}(t)$ for $t \in[0,T]$.

		      Let us next focus on the electric remainder $\tilde{R}^u_{n,\delta,\epsilon}$.
		      We expand the time derivative of the commutator $\pt\cD^u_n = \pt \bigl( \pt^{m-1}(\alpha_n \pt u_n) - \alpha_n \pt^m u_n + \pt^{m-1}(\sigma_n u_n) - \sigma_n \pt^{m-1} u_n \bigr)$.
		      Observe that
		      \begin{align*}
			      \pt \bigl( \pt^{m-1}(\alpha_n \pt u_n) - \alpha_n \pt^m u_n\bigr) & = \pt \left(\sum_{j=1}^{m-1}\binom{m-1}{j}\pt^j \alpha_n \pt^{m-j}u_n\right)                     \\
			                                                                        & = \sum_{j=2}^{m-1} \left(\binom{m-1}{j-1} + \binom{m-1}{j}\right) \pt^j \alpha_n \pt^{m-j+1} u_n \\*
			                                                                        & \seq + \binom{m-1}{m-1} \pt^m \alpha_n \pt u_n+ \binom{m-1}{1} \pt \alpha_n \pt^m u_n            \\
			                                                                        & = \sum_{j=1}^{m}\binom{m}{j} \pt^j \alpha_n \pt^{m-j+1} u_n - (\pt \alpha_n) \pt^m u_n.
		      \end{align*}
		      A similar calculation for $\pt\bigl(\pt^{m-1}(\sigma_n u_n) - \sigma_n \pt^{m-1} u_n\bigr)$ allows us to rewrite $\pt\cD^u_n = \cL^u_n + (m-1)(\pt \alpha_n)U_n - (\pt \sigma_n)\pt^{m-1}u_n$
		      with
		      \[
			      \cL^u_n \coloneqq \sum_{j=2}^m \binom{m}{j} \pt^j \alpha_n \pt^{m-j+1} u_n + \sum_{j=1}^m \binom{m}{j} \pt^j \sigma_n \pt^{m-j}u_n,
		      \]
		      Thus, all remaining terms containing $U_n$ in \eqref{eq:proof-remainder-terms-unshifted}, namely $\sigma_n U_{n,\delta} - \pt(\sigma_n \pt^{m-1}u)_\delta - \cM^n_\delta + m(\pt \alpha_n)U_{n,\delta}$ cancel in the limit with $(m-1)(\pt \alpha_n)U_n - (\pt \sigma_n)\pt^{m-1}u_n$, leaving the lower order term $\cL^u_n$.
		      Thus, $\tilde{R}^u_{n,\delta,\epsilon}$ converges to $\tilde{R}^u_{n,\epsilon}$ in $L^2([0,T]\times\Omega)$ as $\delta \to 0$ given by
		      \begin{align*}
			      \tilde{R}^u_{n,\epsilon} & \coloneqq - (\alpha_n - \alpha)\pt U_\epsilon - \sigma_n U_\epsilon - m (\pt \alpha_n)U_\epsilon             \\
			                               & \seq + \pt(\sigma \pt^{m-1} u)_\epsilon + (\pt^m f)_\epsilon - \pt^m f_n                                     \\
			                               & \seq + \pt (\cD^u)_\epsilon - \cL^u_{n}+ \cM_\epsilon                                                        \\
			                               & \seq - \alpha_n \pt E^u_{n,\epsilon} + \curl E^v_{n,\epsilon} - (\sigma_n + m\pt \alpha_n) E^u_{n,\epsilon},
		      \end{align*}
		      cf. Equations~\eqref{eq:proof-remainder-terms-unshifted} and \eqref{eq:proof-remainder-terms-shifted}.
		      The magnetic terms are treated analogously.
		      Hence, \eqref{eq:proof:gronwall-approx} reduces to
		      \begin{equation}
			      \label{eq:Gronwall1}
			      \tilde{e}^D_{n,\epsilon}(t) \leq \e^{C_G t} \left( \tilde{e}^D_{n,\epsilon}(0) + \int_0^t \Big( \smallnorm{\tilde{R}^u_{n,\epsilon}(s)}_{L^2(\Omega)}^2 + \smallnorm{\tilde{R}^v_{n,\epsilon}(s)}_{L^2(\Omega)}^2 \Big) \dd s \right)
		      \end{equation}
		      for $\delta \to 0$.

		\item Next, we show that $\tilde{e}^D_{n,\epsilon}(0)$ vanishes for $n \to \infty$ and $\epsilon \to 0$. We have
		      \begin{align*}
			      \bignorm{\tilde{D}^u_{n,\epsilon}(0)}_{L^2(\Omega)} & \leq \norm{D^u_{n,\epsilon}(0)}_{L^2(\Omega)} + \norm{E^u_{n,\epsilon}(0)}_{L^2(\Omega)}                                                   \\
			                                                          & \leq \norm{\pt^m u_n(0) - \pt^m u(0)}_{L^2(\Omega)} + \norm{\pt^m u(0) - (\pt^m u)_\epsilon(0)}_{L^2(\Omega)}                              \\
			                                                          & \seq + \!\norm{\pt^m W^u_n(0) \!-\! \pt^m W^u(0)}_{L^2(\Omega)} \!\!+\! \norm{\pt^m W^u(0) \!-\! (\pt^m W^u)_\epsilon(0)}_{L^2(\Omega)}\!.
		      \end{align*}
		      We treat these limits separately.

		      To control the first term, we evaluate the time derivatives of the solutions at $t=0$ recursively by applying the Maxwell system. So, we derive
		      \[
			      \pt^m (u_{0}, v_{0})= ((u_{0})_{(m)},(v_{0})_{(m)}) \quad \text{and} \quad \pt^m (u_{0,n}, v_{0,n})= ((u_{0,n})_{(m)},(v_{0,n})_{(m)}),
		      \]
		      see Definition~\ref{def:CC}, which only depend on the initial data and its spatial derivatives at $t=0$. By \ref{hypothesis-approx}, the initial data converge $(u_{0,n}, v_{0,n}) \to (u_0, v_0)$ in $H^m_\Gamma(\Omega)$, and similarly for the inhomogeneities, coefficients, boundary, and interface data.
		      We thus obtain the convergence of $\pt^m (u_n, v_n)(0) \to \pt^m (u,v)(0)$ in $L^2(\Omega)$.

		      For the second term, we proceed as above, using the extensions to $[-\epsilon_0,T+\epsilon_0]$. 
		      By standard properties of mollifiers, the convergence is uniform in time.
		      Thus,
		      \[
			      \lim_{\epsilon \to 0} \norm{\pt^m u(0) - (\pt^m u)_\epsilon(0)}_{L^2(\Omega)} = 0, \quad \text{and} \quad \lim_{\epsilon \to 0} \norm{\pt^m v(0) - (\pt^m v)_\epsilon(0)}_{L^2(\Omega)} = 0.
		      \]
		      Moreover, since the extensions are chosen such that $W^u_n \to W^u$ in $\cG^{m+1}_\Gamma$, we in particular have $\pt^m W^u_n(0) \to \pt^m W^u(0)$.
		      Likewise, one can show the convergence of the magnetic difference $D^v_{n,\epsilon}(0)$. Hence,
		      \[
			      \lim_{\epsilon \to 0} \limsup_{n \to \infty} \tilde{e}^D_{n,\epsilon}(0) = 0.
		      \]

		\item For fixed $\epsilon \in (0, \epsilon_0)$, we now take the limit $n \to \infty$. Note that $(\alpha_n - \alpha)\pt U_\epsilon \to 0$ in $L^2([0,T] \times \Omega)$.
		      The boundedness of $(u_n, v_n)$ in $\cG_\Gamma^m$, convergence $(u_n, v_n) \to (u, v)$ in $\cG_\Gamma^{m-1}$ from the a-priori estimate \eqref{eq:lwp-apriori} and the Sobolev embeddings shows that the commutator
		      \[
			      \cL^u_{n} \to \cL^u \coloneqq \sum_{j=2}^m \binom{m}{j} \pt^j \alpha \pt^{m-j+1} u + \sum_{j=1}^m \binom{m}{j} \pt^j \sigma \pt^{m-j}u
		      \]
		      converges in $L^2([0,T] \times \Omega)$ for $n \to \infty$.
		      Additionally, since $W^u_n \to W^u$, we have $E^u_{n,\epsilon} \to E^u_\epsilon \coloneqq \pt^m W^u - (\pt^m W^u)_\epsilon$ in $H^1(\Omega)$, and its time derivative in $L^2(\Omega)$.

		      Thus, $\tilde{R}^u_{n,\epsilon}$ tends to
		      \[
			      \tilde{R}^u_{\epsilon} \coloneqq R^u_\epsilon - \alpha \pt E^u_{\epsilon} + \curl E^v_{\epsilon} - (\sigma + m\pt \alpha) E^u_{\epsilon}
		      \]
		      where
		      \begin{equation}\label{eq:proof-defR^u_epsilon}
			      R^u_{\epsilon} \coloneqq \pt(\sigma \pt^{m-1} u)_\epsilon - \sigma U_\epsilon + (\pt^m f)_\epsilon - \pt^m f + \pt (\cD^u)_\epsilon - m (\pt \alpha)U_\epsilon - \cL^u + \cM_\epsilon.
		      \end{equation}
		      in $L^2([0,T] \times \Omega)$ for $n \to \infty$.

		      Finally, we use that $(\pt^m f)_\epsilon \to \pt^m f$ in $L^2([0,T] \times \Omega)$ for $\epsilon \to 0$, which thus vanishes in the limit in \eqref{eq:proof-defR^u_epsilon}.
		      Lemma~\ref{lem:commutator-mollifier} then yields
		      \begin{align*}
			      \lim_{\epsilon \to 0} R^u_{\epsilon} & = (\pt \sigma)\pt^{m-1} u + \pt \cD^u - m (\pt \alpha)U - \cL^u + (\pt \alpha) U                                                    \\
			                                           & = \pt \cD^u - \left(\sum_{j=1}^m \binom{m}{j} \pt^j \alpha \pt^{m-j+1} u + \sum_{j=1}^m \binom{m}{j} \pt^j \sigma \pt^{m-j}u\right) \\
			                                           & \seq + (\pt \alpha) U + (\pt \sigma) \pt^{m-1} u
		      \end{align*}
		      in $L^2([0,T] \times \Omega)$.

		      The summand $\pt (\cD^u)$ is expanded similar to Step~\ref{step:proof-time-convergence-delta-limit}.
		      We rewrite the first term as
		      \begin{align*}
			      \pt \cD^u & = \pt \big( \pt^{m-1}(\alpha \pt u) - \alpha \pt^m u + \pt^{m-1}(\sigma u) - \sigma \pt^{m-1} u \big)                                                                \\
			                & = \left(\sum_{j=1}^m \binom{m}{j} \pt^j \alpha \pt^{m-j+1} u + \sum_{j=1}^m \binom{m}{j} \pt^j \sigma \pt^{m-j}u\right) - (\pt \alpha) U - (\pt \sigma) \pt^{m-1} u,
		      \end{align*}
		      and thus $\lim_{\epsilon \to 0} \int_0^t \smallnorm{R^u_{\epsilon}}_{L^2(\Omega)}^2 \dd \tau = 0$.
		      Furthermore, since $E^u_\epsilon \to 0$ in $H^1(\Omega)$ and $\pt E^u_\epsilon \to 0$ in $L^2(\Omega)$, the additional extension terms in $\tilde{R}^u_\epsilon$ vanish as well.
		      Therefore, $\lim_{\epsilon \to 0} \int_0^t \smallnorm{\tilde{R}^u_{\epsilon}}_{L^2(\Omega)}^2 \dd \tau = 0$.
		      Likewise, one shows the analogous estimate for the magnetic remainder terms $\lim_{\epsilon \to 0} \int_0^t \smallnorm{\tilde{R}^v_{\epsilon}}_{L^2(\Omega)}^2 \dd \tau = 0$.

		\item Therefore, the right-hand side of \eqref{eq:Gronwall1} vanishes and the claim follows, since
		      \begin{align*}
			      \norm{\pt^m u_n(t) - \pt^m u(t)}_{L^2(\Omega)} & \leq \bignorm{\tilde{D}^u_{n,\epsilon}(t)}_{L^2(\Omega)} + \norm{E^u_{n,\epsilon}(t)}_{L^2(\Omega)} \\
			                                                     & \seq + \norm{(\pt^m u)_\epsilon(t) - \pt^m u(t)}_{L^2(\Omega)} \to 0,
		      \end{align*}
		      and analogously for $v$. \qedhere
	\end{enumerate}
\end{proof}

Next, we have to transform the convergence in time into spatial convergence using the div-curl estimate from Theorem~\ref{thm:div-curl-perfect}.

\begin{lemma}[Convergence of Spatial Derivatives]\label{lem:spatial_convergence}
	Let \ref{hypothesis-approx} hold. Suppose that $(u_n,v_n) \to (u,v)$ in $\cG_\Gamma^{m-1}$ and for some $j \in \{1, \dots, m\}$, we have
	\begin{equation}\label{eq:spatial_convergence-assumption}
		\lim_{n \to \infty} \sup_{t \in [0,T]} \Big( \bignorm{\pt^{j} (u_n - u)(t)}_{H_\Gamma^{m-j}} + \bignorm{\pt^{j} (v_n - v)(t)}_{H_\Gamma^{m-j}} \Big) = 0.
	\end{equation}
	Then, the $(j-1)$-th time derivative converges in $H^{m-j+1}_\Gamma$, i.e.,
	\[
		\lim_{n \to \infty} \sup_{t \in [0,T]} \Big( \bignorm{\pt^{j-1} (u_n - u)(t)}_{H_\Gamma^{m-j+1}} + \bignorm{\pt^{j-1} (v_n - v)(t)}_{H_\Gamma^{m-j+1}} \Big) = 0.
	\]
\end{lemma}

\begin{proof}
	Let $r \coloneqq m - j + 1$. We define the differences $D^u_n \coloneqq \pt^{j-1} u_n - \pt^{j-1} u$ and $D^v_n \coloneqq \pt^{j-1} v_n - \pt^{j-1} v$. We have to show that
	\[
		e^D_n(t) \coloneqq \norm{D^u_n(t)}_{H_\Gamma^r} + \norm{D^v_n(t)}_{H_\Gamma^r} \to 0
	\]
	uniformly for $t \in [0,T]$ and $n \to \infty$.

	In order to apply the div-curl Theorem~\ref{thm:div-curl-perfect}, the fields must belong to $H_{t0}(\curl, \Omega)$, respectively $H(\curl, \Omega)$. Since $D^u_n$ and $D^v_n$ satisfy inhomogeneous boundary and interface conditions, we again extend these errors to the interior.

	As in the proof of the previous lemma there exist extensions $(W^u_n, W^v_n)$ of $(\bndInh_n, \kappa_n, \varkappa_n)$ and $(W^u, W^v)$ of $(\bndInh, \kappa, \varkappa)$ in $\cG^{m+1}_\Gamma$ converging in $\cG^{m+1}_\Gamma$ following the approximation in \ref{def:approx-data}.
	Let $E^u_n \coloneqq \pt^{j-1}(W^u_n - W^u)$ and $E^v_n \coloneqq \pt^{j-1}(W^v_n - W^v)$. Since the extensions belong to $\cG^{m+1}_\Gamma$, these $(j-1)$-th time derivatives belong to $H^{(m+1)-(j-1)}_\Gamma(\Omega) = H^{r+1}_\Gamma(\Omega)$.
	We have
	\[
		\sup_{t \in [0,T]} \Big( \norm{E^u_n(t)}_{H^{r+1}_\Gamma} + \norm{E^v_n(t)}_{H^{r+1}_\Gamma} \Big) \to 0
	\]
	as $n \to \infty$.

	We define the shifted differences $\tilde{D}^u_n \coloneqq D^u_n - E^u_n$ and $\tilde{D}^v_n \coloneqq D^v_n - E^v_n$. By construction, these satisfy exactly homogeneous tangential traces and jumps:
	\[
		\tilde{D}^u_n \times \nu = 0 \text{ on } \partial\Omega, \qquad \llbracket \tilde{D}^u_n \times \nu^\Gamma \rrbracket = 0, \qquad \llbracket \tilde{D}^v_n \times \nu^\Gamma \rrbracket = 0.
	\]

	We start by deriving the Maxwell equations satisfied by the unshifted differences $(D^u_n, D^v_n)$.
	For $(u_n, v_n)$ we have
	\begin{alignat*}{-1}
		\alpha_n \pt^j u_n & =   &  & \curl \pt^{j-1} v_n &  & - \sigma_n \pt^{j-1} u_n &  & - \pt^{j-1} f_n &  & - \cL^u_{n}, \\
		\beta_n \pt^j v_n  & = - &  & \curl \pt^{j-1} u_n &  & - \rho_n \pt^{j-1} v_n   &  & - \pt^{j-1} g_n &  & - \cL^v_{n},
	\end{alignat*}
	with $\cL^u_{n} \coloneqq [\pt^{j-1}, \alpha_n \pt] u_n + [\pt^{j-1}, \sigma_n] u_n$ and $\cL^v_{n} \coloneqq [\pt^{j-1}, \beta_n \pt] v_n+ [\pt^{j-1}, \rho_n] v_n$.

	Similarly, the limit $(u,v)$ satisfies
	\begin{alignat*}{-1}
		\alpha \pt^j u & =   &  & \curl \pt^{j-1} v &  & - \sigma \pt^{j-1} u &  & - \pt^{j-1} f - \cL^u, \\
		\beta \pt^j v  & = - &  & \curl \pt^{j-1} u &  & - \rho \pt^{j-1} v   &  & - \pt^{j-1} g - \cL^v,
	\end{alignat*}
	with commutators $\cL^u \coloneqq [\pt^{j-1}, \alpha \pt] u + [\pt^{j-1}, \sigma] u$ and $\cL^v \coloneqq [\pt^{j-1}, \beta \pt] v+ [\pt^{j-1}, \rho] v$.
	Thus, their differences satisfy
	\begin{alignat*}{-1}
		\pt(\alpha_n D^u_n) & =   &  & \curl D^v_n &  & - J^u_n, \\
		\pt(\beta_n D^v_n)  & = - &  & \curl D^u_n &  & - J^v_n,
	\end{alignat*}
	where we set
	\begin{alignat*}{-1}
		J^u_n & \!\coloneqq\! (\sigma_n - \pt \alpha_n)D^u_n &  & + (\alpha_n - \alpha)\pt^j u &  & + (\sigma_n - \sigma)\pt^{j-1} u &  & + \pt^{j-1}(f_n - f) &  & + (\cL^u_{n} - \cL^u), \\*
		J^v_n & \!\coloneqq\! (\rho_n - \pt\beta_n) D^v_n    &  & + (\beta_n - \beta)\pt^j v   &  & + (\rho_n - \rho)\pt^{j-1} v     &  & + \pt^{j-1}(g_n - g) &  & + (\cL^v_{n} - \cL^v).
	\end{alignat*}

	Substituting the shifted variables, the system becomes
	\begin{alignat}{-1}
		\label{eq:diff_sys_u}\pt(\alpha_n \tilde{D}^u_n) & =   &  & \curl \tilde{D}^v_n &  & - \tilde{J}^u_n, \\
		\label{eq:diff_sys_v}\pt(\beta_n \tilde{D}^v_n)  & = - &  & \curl \tilde{D}^u_n &  & - \tilde{J}^v_n,
	\end{alignat}
	where the modified remainders are given by
	\begin{alignat*}{-1}
		\tilde{J}^u_n & \coloneqq J^u_n &  & + \pt(\alpha_n E^u_n) &  & - \curl E^v_n, \\
		\tilde{J}^v_n & \coloneqq J^v_n &  & + \pt(\beta_n E^v_n)  &  & + \curl E^u_n.
	\end{alignat*}

	Because $\tilde{D}^u_n$ and $\tilde{D}^v_n$ possess strictly zero tangential jumps, their curls have no surface divergence. Due to Lemma~8.1 in \cite{schnaubelt2022}, the normal jumps are thus given by
	\begin{alignat*}{-1}
		\llbracket \alpha_n \tilde{D}^u_n \cdot \nu^\Gamma \rrbracket(t) & = \llbracket \alpha_n \tilde{D}^u_n \cdot \nu^\Gamma \rrbracket(0) &  & - \int_0^t \llbracket \tilde{J}^u_n(\tau) \cdot \nu^\Gamma \rrbracket \dd \tau, \\
		\llbracket \beta_n \tilde{D}^v_n \cdot \nu^\Gamma \rrbracket(t)  & = \llbracket \beta_n \tilde{D}^v_n \cdot \nu^\Gamma \rrbracket(0)  &  & - \int_0^t \llbracket \tilde{J}^v_n(\tau) \cdot \nu^\Gamma \rrbracket \dd \tau,
		\intertext{and}
		\tr_n (\beta_n \tilde{D}^v_n)(t)                                 & = \tr_n (\beta_n \tilde{D}^v_n)(0)                                 &  & - \int_0^t \tr_n \tilde{J}^v_n(\tau) \dd \tau.
	\end{alignat*}
	Thus, in view of Theorem~\ref{thm:div-curl-perfect} (corresponding to $v$ there), we define the correction terms
	\begin{equation}\label{eq:def-tilde-C^uv_n}
		\tilde{C}^u_n \coloneqq -\alpha_n(0) \tilde{D}^u_n(0) + \int_0^t \tilde{J}^u_n(\tau)\dd \tau, \quad
		\tilde{C}^v_n \coloneqq -\beta_n(0) \tilde{D}^v_n(0) + \int_0^t \tilde{J}^v_n(\tau)\dd \tau.
	\end{equation}
	Proposition~2.2.30 in \cite{ciarlet18}, see Lemma~\ref{lem:characterization-H-div-H-curl}, then shows that $\alpha_n \tilde{D}^u_n + \tilde{C}^u_n \in H(\div,\Omega)$ and $\beta_n \tilde{D}^v_n + \tilde{C}^v_n \in H_{n0}(\div,\Omega)$.
	We further obtain $\tilde{D}_n^u \in H_{t0}(\curl, \Omega)$ and $\tilde{D}_n^v \in H(\curl, \Omega)$.

	Theorem~\ref{thm:div-curl-perfect} then yields for the shifted energy $\tilde{e}^D_n(t) \coloneqq \bignorm{\tilde{D}^u_n(t)}_{H_\Gamma^r} + \bignorm{\tilde{D}^v_n(t)}_{H_\Gamma^r}$
	\begin{align}\label{eq:div-curl-expansion}
		\tilde{e}^D_n(t) & \lesssim \bignorm{\tilde{D}^u_n(t)}_{L^2(\Omega)} + \bignorm{\curl \tilde{D}^u_n(t)}_{H_\Gamma^{r-1}} + \bignorm{\div(\alpha_n \tilde{D}^u_n(t))}_{H_\Gamma^{r-1}} + \bignorm{\tilde{C}^u_n(t)}_{H_\Gamma^r} \nonumber \\
		                 & \quad + \bignorm{\tilde{D}^v_n(t)}_{L^2(\Omega)} + \bignorm{\curl \tilde{D}^v_n(t)}_{H_\Gamma^{r-1}} + \bignorm{\div(\beta_n \tilde{D}^v_n(t))}_{H_\Gamma^{r-1}} + \bignorm{\tilde{C}^v_n (t)}_{H_\Gamma^r} .
	\end{align}
	To show convergence, we now consider the terms on the right-hand side individually. Again, we only consider $u$ as the estimates for $v$ follow in the same way.
	\begin{enumerate}[label = \arabic*), wide, labelwidth=0pt, labelindent=0pt, itemsep=0.5em, topsep=0.5em]
		\item First, the $L^2$ norms converge. We have $\bignorm{\tilde{D}^u_n(t)}_{L^2(\Omega)} \leq \norm{D^u_n(t)}_{L^2(\Omega)} + \norm{E^u_n(t)}_{L^2(\Omega)}$. By assumption $(u_n,v_n) \to (u,v)$ in $\cG_\Gamma^{m-1}$.
		      Thus, in particular $D^u_n(t) \to 0$ in $L^2(\Omega)$. Since $E^u_n(t) \to 0$ in $H^{r+1}_\Gamma$, we obtain $\bignorm{\tilde{D}^u_n(t)}_{L^2(\Omega)} \to 0$ uniformly for $t \in [0,T]$.

		\item \label{step:proof-spatial-2} Next, we consider the correction terms $(\tilde{C}^u_n, \tilde{C}^v_n)$.
		      We decompose $\tilde{J}^u_n$ as follows. To bound the commutator difference $\cL^u_n - \cL^u$, we extract its highest-order term $(j-1)(\pt \alpha_n)\pt^{j-1} u_n - (j-1)(\pt \alpha)\pt^{j-1} u$.
		      Using $D^u_n = \tilde{D}^u_n + E^u_n$, we move the $(j-1)(\pt \alpha_n)\tilde{D}^u_n$ component into the Gronwall term. \footnote{This highest-order term was overlooked in a draft of the paper. The AI tool \textit{Claude Fable 5.1} caught that mistake.}
		      Thus, we write
		      \[
			      \tilde{J}^u_n = (\sigma_n + (j-2)\pt \alpha_n)\tilde{D}^u_n + \tilde{F}^u_n,
		      \]
		      with
		      \begin{align*}
			      \tilde{F}^u_n \coloneqq & (\alpha_n - \alpha)\pt^j u + (\sigma_n - \sigma)\pt^{j-1} u + \pt^{j-1}(f_n - f) \\
			                              & + \big(\cL^u_{n} - \cL^u - (j-1)(\pt \alpha_n)\tilde{D}^u_n\big)                 \\
			                              & + \sigma_n E^u_n + \alpha_n \pt E^u_n - \curl E^v_n,
		      \end{align*}
		      where we have subtracted $(j-1)(\pt \alpha_n)\tilde{D}^u_n$ such that the commutator difference
		      \begin{align*}
			      \cL^u_n - \cL^u - (j-1)(\pt \alpha_n)\tilde{D}^u_n
			       & = (j-1)(\pt \alpha_n) E^u_n + (j-1)(\pt \alpha_n - \pt \alpha)\pt^{j-1} u                                                            \\
			       & \quad + \sum_{\ell=2}^{j-1} \binom{j-1}{\ell} \big[ (\pt^\ell \alpha_n)\pt^{j-\ell} u_n - (\pt^\ell \alpha)\pt^{j-\ell} u \big]      \\
			       & \quad + \sum_{\ell=1}^{j-1} \binom{j-1}{\ell} \big[ (\pt^\ell \sigma_n)\pt^{j-1-\ell} u_n - (\pt^\ell \sigma)\pt^{j-1-\ell} u \big].
		      \end{align*}
		      no longer depends on the difference $D^u_n = \pt^{j-1}(u_n - u)$ of order $j-1$ in time.
		      The first term in $\tilde{J}^u_n$ will be treated by a Gronwall argument.
		      To bound the norm of the integral of $\tilde{F}^u_n$ in $H_\Gamma^r$, we note that $\pt^j u \in H_\Gamma^{r-1}$, which prevents a direct bound. However, integrating by parts in time yields
		      \begin{equation}\label{eq:proof:ibp}
			      \int_0^t (\alpha_n - \alpha)\pt (\pt^{j-1} u) \dd \tau = [(\alpha_n - \alpha)\pt^{j-1} u]_0^t - \int_0^t \pt(\alpha_n - \alpha)\pt^{j-1} u \dd \tau.
		      \end{equation}
		      Since $\pt^{j-1} u \in C([0,T]; H_\Gamma^r)$ and $\pt\alpha_n \to \pt\alpha$ in $\cC^{\hat m}$, by Lemma~2.1 in \cite{schnaubelt2022}, the $H_\Gamma^r$-norm of the right-hand side converges to zero uniformly.
		      The remaining terms in $\tilde{F}^u_n$, except the inhomogeneity and the extension terms, similarly vanish after integration in $H_\Gamma^r$ due to the convergence of the coefficients and commutators, since $u_n\to u$ in $\cG_\Gamma^{m-1}$ and $u \in \cG^{m}_\Gamma$.
		      For the inhomogeneity, since $f \in H^m([0,T]\times\Omega)$, we have $\pt^{j-1}(f_n - f) \in L^2([0,T]; H^r(\Omega))$ and
		      \[
			      \int_0^t \norm{\pt^{j-1}(f_n - f)(\tau)}_{H^r_\Gamma} \dd \tau \leq \sqrt{T} \norm{\pt^{j-1}(f_n - f)}_{L^2([0,T]; H^r_\Gamma)} \to 0
		      \]
		      as $n \to \infty$. Furthermore, the terms $\sigma_n E^u_n, \alpha_n \pt E^u_n, \curl E^v_n$ vanish uniformly in $H^r_\Gamma$ because $E^u_n, E^v_n \to 0$ in $H^{r+1}_\Gamma$ and $\pt E^u_n \to 0$ in $H^r_\Gamma$.
		      Thus, keeping the norm strictly outside the time integral for $\tilde{F}^u_n$, we obtain
		      \begin{align*}
			      \bignorm{\tilde{C}^u_n(t)}_{H_\Gamma^r} & \lesssim \bignorm{\alpha_n(0) \tilde{D}^u_n(0)}_{H_\Gamma^r} + \biggnorm{\int_0^t \tilde{F}^u_n(\tau) \dd \tau}_{H_\Gamma^r} + \int_0^t \bignorm{\tilde{D}^u_n(\tau)}_{H_\Gamma^r} \dd \tau \\
			                                              & \leq \tilde{\delta}_{u,n} + \int_0^t \tilde{e}^D_n(\tau) \dd \tau,
		      \end{align*}
		      for $\tilde{\delta}_{u,n} \coloneqq (1+\norm{\alpha_n(0)}_{H_\Gamma^{\tilde m}})\bignorm{\tilde{D}^u_n(0)}_{H_\Gamma^r} + \sup_{t \in [0,T]} \bignorm{\int_0^t \tilde{F}^u_n(\tau) \dd \tau}_{H_\Gamma^r} \to 0$ as $n \to \infty$.

		\item We consider the divergence terms next.
		      Integrating \eqref{eq:diff_sys_u} and substituting \eqref{eq:def-tilde-C^uv_n} yields $\alpha_n \tilde{D}^u_n = - \tilde C_n^u -\curl \int_0^t \tilde{D}^u_n\dd \tau$ on $\Omega_\pm$. Therefore, taking the norm of the divergence yields
		      \[
			      \bignorm{\div(\alpha_n \tilde{D}^u_n(t))}_{H_\Gamma^{r-1}} \lesssim \bignorm{\tilde{C}^u_n(t)}_{H_\Gamma^r}.
		      \]
		      This is the term we have already considered in the previous step. Thus,
		      \[
			      \bignorm{\div(\alpha_n \tilde{D}^u_n(t))}_{H_\Gamma^{r-1}} \lesssim \tilde{\delta}_{u,n} + \int_0^t \tilde{e}^D_n(\tau) \dd \tau.
		      \]
		      In the same way, we bound $\bignorm{\div(\beta_n \tilde{D}^v_n)}_{H_\Gamma^{r-1}}$.

		\item Finally, let us consider the curl terms.
		      By \eqref{eq:diff_sys_v}, and the definition of $\tilde{J}^v_n$, we write as in Step~\ref{step:proof-spatial-2}
		      \[
			      \curl \tilde{D}^u_n = -\pt(\beta_n \tilde{D}^v_n) - \tilde{J}^v_n = -\beta_n \pt \tilde{D}^v_n - (\rho_n + (j-1)\pt \beta_n)\tilde{D}^v_n - \tilde{F}^v_n,
		      \]
		      with
		      \begin{align*}
			      \tilde{F}^v_n \coloneqq & (\beta_n - \beta)\pt^j v + (\rho_n - \rho)\pt^{j-1} v + \pt^{j-1}(g_n - g) \\
			                              & + \big(\cL^v_{n} - \cL^v - (j-1)(\pt \beta_n)\tilde{D}^v_n\big)            \\
			                              & + \rho_n E^v_n + \beta_n \pt E^v_n + \curl E^u_n.
		      \end{align*}

		      This allows us to trade space against time regularity of $\pt \tilde{D}^v_n = \pt^j (v_n - v) - \pt E^v_n$.
		      By assumption~\eqref{eq:spatial_convergence-assumption}, $\pt^j (v_n - v) \to 0$ uniformly in $H_\Gamma^{m-j} = H_\Gamma^{r -1}$ for $n \to \infty$.
		      Since $\pt E^v_n \to 0$ in $H^r_\Gamma \subset H^{r-1}_\Gamma$, we deduce $\pt \tilde{D}^v_n \to 0$ in $H^{r-1}_\Gamma$.
		      Additionally, $\pt^j v \in H_\Gamma^{m-j}$, so the uniform convergence of $\beta_n$, $\rho_n$, the inhomogeneities, the extensions, and the commutators implies that $\bignorm{\tilde{F}^v_n}_{H_\Gamma^{m-j}} \to 0$. Thus,
		      \[
			      \sup_{t \in [0,T]} \bignorm{\curl \tilde{D}^u_n(t)}_{H_\Gamma^{r-1}} \to 0 \quad \text{as } n \to \infty.
		      \]
		      Based on \eqref{eq:diff_sys_u}, an analogous argument shows $\sup_{t \in [0,T]} \bignorm{\curl \tilde{D}^v_n(t)}_{H_\Gamma^{r-1}} \to 0$.
	\end{enumerate}

	Inserting the above estimates into \eqref{eq:div-curl-expansion} we arrive at the inequality
	\[
		\tilde{e}^D_n(t) \leq \Delta_n + C \int_0^t \tilde{e}^D_n(\tau) \dd \tau,
	\]
	where $C > 0$ and
	\[
		\Delta_n \coloneqq \sup_{t \in [0,T]} \left(\Bignorm{(\tilde{D}^u_n,\tilde{D}^v_n)(t)}_{L^2(\Omega)} + \Bignorm{\curl(\tilde{D}^u_n,\tilde{D}^v_n)(t)}_{H^{r-1}_\Gamma}\right) + 2 (\tilde{\delta}_{u,n} + \tilde{\delta}_{v,n}).
	\]
	Since $\Delta_n \to 0$ as $n \to \infty$, Gronwall's inequality yields
	\[
		\sup_{t \in [0,T]} \tilde{e}^D_n(t) \leq \Delta_n\, \mathrm{e}^{CT} \to 0 \quad \text{for } n \to \infty.
	\]
	Finally, for the unshifted variables, we have
	\[
		\sup_{t \in [0,T]} e^D_n(t) \leq \sup_{t \in [0,T]} \Big(\tilde{e}^D_n(t) + \norm{E^u_n(t)}_{H^r_\Gamma} + \norm{E^v_n(t)}_{H^r_\Gamma}\Big) \to 0 \quad \text{for } n \to \infty,
	\]
	and the claim follows.
\end{proof}

Now, we combine Lemma~\ref{lem:time-convergence} and Lemma~\ref{lem:spatial_convergence} to obtain our approximation result.

\begin{theorem}\label{thm:approx}
	Let $m \in \N_0, k,a \in \N$ with $\tilde m + a + 2 \leq k$, and $\partial\Omega \subseteq \partial\omega \in C^{k+1}$.
	Furthermore, let $\alpha, \beta \in \cC^{\tilde m}_\eta, \pt\alpha, \pt\beta \in \cC^{\hat m}, \sigma, \rho \in \cC^{\tilde m}, f,g \in H^m([0,T] \times \Omega_\pm)$, $\bndInh \in T^m(\partial \Omega)$, $\kappa,\varkappa \in J^m(\Gamma)$, and $u_0,v_0 \in H_\Gamma^m(\Omega)$ satisfy the compatibility conditions $(CC_{\alpha,\beta, \sigma, \rho, f, g, \bndInh, \kappa, \varkappa, m})$ of order~$m$.
	By Theorem 3.1 in \cite{schnaubelt2022}, see Proposition~\ref{prop:lwp-apriori}, there exists a solution $(u,v) \in \cG_\Gamma^m([0,T])$ to \eqref{eq:maxwell-lem}--\eqref{eq:maxwell-lem4}.

	Choose $\alpha_n, \beta_n, \sigma_n, \rho_n \in C^\infty([0,T] \times \Omega_\pm)$, $\alpha_n, \beta_n$ uniformly positive definite, $f_n, g_n \in C^\infty([0,T] \times \Omega_\pm)$, $\bndInh_n \in C^\infty([0,T] \times \partial\Omega)$, $\kappa_n, \varkappa_n \in C^\infty([0,T] \times \Gamma)$ satisfying \ref{def:approx-data}. Then there exist $u_{0,n}, v_{0,n} \in H_\Gamma^{m+a}(\Omega)$ satisfying the compatibility conditions $(CC_{\alpha_n,\beta_n, \sigma_n, \rho_n, f_n, g_n, \bndInh_n, \kappa_n, \varkappa_n, m+a})$ of order~$m+a$,
	such that there exists a corresponding solution $(u_n,v_n)\in \cG_\Gamma^{m+a}([0,T])$ of \eqref{eq:maxwell-lem} with
	\[
		(u_n,v_n) \to (u,v) \qquad \text{in }\,\cG_\Gamma^m([0,T])
	\]
	for $n \to \infty$.
\end{theorem}

\begin{proof}
	By Remark~\ref{rem:approx-m-geq-1} we may assume $m \geq 1$ in the following.
	Given $\alpha,\beta, \sigma, \rho, f, g$, boundary and jump data $\bndInh, \kappa, \varkappa$, and initial values $u_0,v_0$, we choose an approximating sequence as outlined in \ref{hypothesis-approx}. Recall that
	\[
		\norm{u}_{\cG_\Gamma^m([0,T] \times \Omega)} \coloneqq \sup_{t \in [0,T]} \sum_{j=0}^m \bignorm{\pt^j u(t)}_{H^{m-j}(\Omega_\pm)}
	\]
	and consider the differences $(u_n - u, v_n - v)$.

	Firstly, we consider the full time derivatives, corresponding to the $j=m$ term in the $\cG_\Gamma^m$-norm.
	Lemma~\ref{lem:time-convergence} yields
	\[
		\lim_{n \to \infty} \sup_{t \in [0,T]} \Big( \norm{\pt^m (u_n - u)(t)}_{L^2(\Omega)} + \norm{\pt^m (v_n - v)(t)}_{L^2(\Omega)} \Big) = 0.
	\]

	Secondly, we use the time convergence as a starting point to iteratively recover the spatial derivatives using Lemma~\ref{lem:spatial_convergence}.
	Considering $j=m$ down to $j=0$, we show that for every $j \in \{0, \dots, m\}$, we have
	\[
		\lim_{n \to \infty} \sup_{t \in [0,T]} \Big( \smallnorm{\pt^j (u_n - u)(t)}_{H^{m-j}(\Omega)} + \smallnorm{\pt^j (v_n - v)(t)}_{H^{m-j}(\Omega)} \Big) = 0,
	\]
	which then shows
	\[
		\lim_{n \to \infty} \Big( \norm{u_n - u}_{\cG_\Gamma^m([0,T] \times \Omega)} + \norm{v_n - v}_{\cG_\Gamma^m([0,T] \times \Omega)} \Big) = 0. \qedhere
	\]
\end{proof}

\begin{remark}\label{rem:quasilinear}
	Theorem~\ref{thm:approx} is an approximation result for the nonautonomous, linear system under the assumption \eqref{eq:additional-reg}. However, we can apply it after differentiating the quasilinear system
	\begin{alignat*}{-1}
		\begin{split}
			\pt \left(\alpha(u,x) u\right) & = \hphantom{-}\curl v - \tilde \sigma(t,x)u - f, \\
			\pt\left(\beta(v,x) v\right)   & = - \curl u - \tilde\rho(t,x)v\hspace{1pt} - g,
		\end{split}
		\hspace{0.5cm}                               &                                                                         & t\geq 0, x \in \Omega_\pm,                                   \\
		u \times \nu(x)                              & = \bndInh,                                                              &                            & t\geq 0, x \in \partial \Omega, \\
		\llbracket u \times \nu^\Gamma(x) \rrbracket & =\kappa,\qquad\llbracket v \times \nu^\Gamma(x) \rrbracket = \varkappa,
		                                             &                                                                         & t\geq 0, x \in \Gamma,                                       \\
		u(0,x)                                       & = u_0(x),  \hspace{8.45mm} v(0,x)= v_0(x),\qquad                        &                            & x \in \Omega_\pm,
	\end{alignat*}
	with coefficients $\alpha(x,u(t,x)), \beta(x,v(t,x))$ depending on $(u,v)$, where usually one only assumes
	\begin{equation}
		\label{eq:quasilinear-condition}
		\partial_{\xi_i} \alpha, \partial_{\xi_i} \beta \in \cC^{\tilde m},
	\end{equation}
	see for example (2.2) in \cite{lasiecka2019}.
	It is a crucial tool in the study of long-time behavior of the quasilinear system, see \cite{nutt2026-2}.

	Given coefficients and initial values satisfying \eqref{eq:maxwell-lem-coeff}, \eqref{eq:maxwell-lem-IVs} and \eqref{eq:quasilinear-condition} for some $m\leq k$ and $r \in \{0,\dots,m\}$, we may consider the $r$-times time differentiated system
	\begin{alignat*}{-1}
		\begin{split}
			\alphad(t,x) \pt \left(\pt^r u(t,x)\right) & = \hphantom{-} \curl \pt^r v(t,x) - \tilde \sigma(t,x)\pt^r u(t,x) \hspace{-0.5pt}- f_r(t,x), \\
			\betad(t,x) \pt\left(\pt^r v(t,x)\right)   & = - \curl \pt^r u(t,x) - \tilde \rho(t,x)\pt^r v(t,x) \hspace{0.5pt}- g_r(t,x),
		\end{split}
		\quad                                           &                                                                                                    & t\geq 0, x \in \Omega_\pm,                                         \\
		\pt^r u(t,x) \times \nu(x)                      & = \pt^r \bndInh(t,x),                                                                              &                            & t\geq 0, x \in \partial \Omega,\notag \\
		\llbracket \pt^r u \times \nu^\Gamma \rrbracket & = \pt^r \kappa(t,x), \quad \llbracket \pt^r v \times \nu^\Gamma \rrbracket = \pt^r \varkappa(t,x), &                            & t\geq 0, x \in \Gamma,
	\end{alignat*}
	where
	\begin{alignat*}{-1}
		\alphad(t,x) & = \mathrlap{\alpha(x,u(t,x)) + \partial_\xi \alpha(x,\xi)\big|_{\xi=u(t,x)} u(t,x) \!= \alpha(x,u) \hspace{-1pt}+ \left(\sum_{\ell=1}^3\partial_{\xi_j}\alpha_{i\ell}(x,u)u_\ell\right)_{\!\!ij} \!\!\!,} &  & \hspace*{\linewidth}                                                                                                                                      \\
		\betad(t,x)  & = \mathrlap{\beta(x,v(t,x)) + \partial_\xi \beta(x,\xi)\big|_{\xi=v(t,x)} v(t,x)= \beta(x,v) + \left(\sum_{\ell=1}^3\partial_{\xi_j}\beta_{i\ell}(x,v)v_\ell\right)_{\!\!ij}\!\!\!,}                                                                                                                                                                                     \\
		f_r(t,x)     & = \pt^r f(t,x)                                                                                                                                                                                            &  & \! +\! \sum_{j=1}^r\! \binom{r}{j} \pt^j \tilde \sigma(t,x) \pt^{r-j} \!u(t,x) \!+\! \sum_{j=1}^{r}\! \binom{r}{j} \pt^j \alphad(t,x) \pt^{r-j+1} u(t,x), \\
		g_r(t,x)     & = \pt^r g(t,x)                                                                                                                                                                                            &  & \!+\! \sum_{j=1}^r\! \binom{r}{j} \pt^j \tilde \rho(t,x) \pt^{r-j} v(t,x) \!+\! \sum_{j=1}^{r}\! \binom{r}{j} \pt^j \betad(t,x) \pt^{r-j+1} v(t,x).
	\end{alignat*}
	Hence, the time derived system has initial data and inhomogeneities which possess regularity of level $m-r$, with coefficients
	\[
		\alphad, \betad \in \cC^{\tilde m} \quad \text{and thus} \quad \pt \alphad, \pt \betad \in \cC^{\tilde m-1}.
	\]
	In this setting we can apply Theorem~\ref{thm:approx}, assuming that $\alphad, \betad$ are uniformly positive definite, to show convergence of approximations to $\pt^r (u,v)$ in $\cG^{m-r}_\Gamma$ for $r \geq 1$ and $m \in \N$.
	Moreover, if $m \leq 2$ then $\cC^{\tilde m-1} = \cC^{\hat m}$ and we obtain convergence in $\cG^{m}_\Gamma$, which we use in \cite{nutt2026-2}.

	In order to obtain convergence in $\cG^{m}_\Gamma$ for general $m\in\N$, one would have to approximate $\alpha(u), \beta(v)$ preserving the quasilinear structure, which is not possible with our current approach, since we assume $\alpha_n,\beta_n$ to be independent of the solution.
\end{remark}

\appendix
\section{Div-Curl estimates and related results}
\label{sec:A1}
\subsection{Properties of the Curl Operator}
We first list basic mapping properties of the gradient, divergence and curl. They can be found in Theorem~2.7 of
\cite{cessenat1996}, Chapter~IX.1.3 of \cite{dautray1990}, or follow similarly. Unless otherwise specified, let $\Omega \subset \R^3$ be an open, bounded, multiply connected set with a Lipschitz boundary.
We partly assume more boundary regularity as in \cite{nutt2026} and \cite{nutt2024}.
However, the regularity assumptions may not be optimal.

\begin{lemma}\label{lem:mapping}
	Let $\Omega\subset \R^3$ be open, bounded and multiply connected with Lipschitz boundary $\partial\Omega$.
	In the following diagram the image of each operator is of finite codimension and contained in the kernel of the
	succeeding operator:
	\begin{align*}
		H^1(\Omega)   & \xrightarrow{\nabla} \makebox[0pt][l]{$H(\curl)$}\hphantom{H_{t0}(\curl)} \xrightarrow{\curl}
		\makebox[0pt][l]{$H(\div)$}\hphantom{H_{n0}(\div)} \xrightarrow{\div} L^2(\Omega),                                  \\
		H_c^1(\Omega) & \xrightarrow{\nabla} H_{t0}(\curl) \xrightarrow{\curl} H_{n0}(\div) \xrightarrow{\div} L^2(\Omega),
	\end{align*}
	where $H_c^1(\Omega)\coloneqq\big\{f\in H^1(\Omega)\bigm|\tr_{\Gamma_k}f\text{ constant for each } k\in\{0,\dots K\} \big\}
		\supset H^1_0(\Omega)$.
\end{lemma}

The following characterization of $H(\div_\alpha)$ and $H(\curl_\alpha)$ in an interface setting can be found in Propositions~2.2.30 and 2.2.32 in \cite{ciarlet18}.

\begin{lemma}\label{lem:characterization-H-div-H-curl}
	Let $\alpha \in L^\infty(\Omega;\R^{3 \times 3})$ invertible with $\alpha^{-1} \in L^\infty(\Omega;\R^{3 \times 3})$ and $\Omega = \upsilon \cupdot \Gamma \cupdot \omega$.
	Then,
	\begin{alignat*}{-1}
		\hspace{-0.8cm}\text{a)}\hspace{0.8cm} &  & H(\div_\alpha , \Omega) = \big\{f \in L^2(\Omega) \,\big\vert\,  & \alpha f|_{\upsilon} \in H(\div,\upsilon), \alpha f|_{\omega} \in H(\div,\omega),       \\
		                                       &  &                                                                  & \llbracket\alpha f \cdot \nu^\Gamma \rrbracket = 0 \text{ in } H^{-1/2}(\Gamma)\big\},  \\
		\hspace{-0.8cm}\text{b)}\hspace{0.8cm} &  & H(\curl_\alpha , \Omega) = \big\{f \in L^2(\Omega) \,\big\vert\, & \alpha f|_{\upsilon} \in H(\curl,\upsilon), \alpha f|_{\omega} \in H(\curl,\omega),     \\
		                                       &  &                                                                  & \llbracket\alpha f \times \nu^\Gamma \rrbracket = 0 \text{ in } H^{-1/2}(\Gamma)\big\}.
	\end{alignat*}
\end{lemma}

We end this section by proving the following estimates for normal traces of the curl operator.

\begin{lemma}\label{Lem:surfaceCurl}
	For $f \in H_\Gamma^1(\Omega)$ we can estimate
	\begin{enumerate}[label =\alph*)]
		\item the normal surface curl by
		      \[
			      \norm{\nu \cdot( \curl f)}_{H^{-1}(\partial\Omega)} \lesssim \norm{\nu \times f}_{L^2(\partial\Omega)},
		      \]
		\item and the normal jump at the interface by
		      \[
			      \norm{\llbracket\nu^\Gamma \cdot( \curl f) \rrbracket}_{H^{-1}(\Gamma)} \lesssim \norm{ \llbracket \nu^\Gamma \times f \rrbracket}_{L^2(\Gamma)}.
		      \]
	\end{enumerate}
\end{lemma}

\begin{proof}
	a) We localize $f$ near $\partial \Omega$ by a cutoff $\chi_{\partial \Omega} \in C^\infty$ and apply Lemma~4.8 from \cite{nutt2024}. It is proven as part b) below. See also Chapter~3.1 of \cite{ciarlet18}.

	b) Likewise, we choose $\chi_\Gamma \in C_c^\infty(\cN_a^\Gamma)$, for a small $a>0$ such that $\cN_a^\Gamma \cap \partial \Omega = \emptyset$, with $0 \leq \chi^\Gamma \leq 1$ and $\chi_\Gamma = 1$ on $\cN_{a/2}^\Gamma$ and essentially repeat the proof in \cite{nutt2024}:

	Let $\phi \in H^1(\Gamma)$.
	This function can be extended in $H^{3/2}(\Omega)$ with $\norm{\phi}_{H^{3/2}(\Omega)} \lesssim \norm{\phi}_{H^1(\Gamma)}$.
	Take $\phi_m \in C^\infty(\overline\Omega)$ with $\phi_m \to \phi$ in $H^{3/2}(\Omega)$ and $f_n \in C^\infty(\overline\Omega \setminus \Gamma)$ with $f_n \to f$ in $H^1_\Gamma(\Omega)$.
	Then the divergence theorem and integration by parts yield
	\begin{align*}
		 & \abs{\oint_\Gamma \llbracket\nu^\Gamma \cdot \curl f_n \rrbracket \phi_m \dd \varsigma}  = \abs{\oint_\Gamma \llbracket\nu^\Gamma \cdot \curl(\chi_\Gamma f_n) \rrbracket \phi_m \dd \varsigma}                                           \\
		 & \quad= \abs{(\curl(\chi_\Gamma f_n)|\nabla\phi_m)_{L^2(\upsilon)} - (\curl(\chi_\Gamma f_n)|\nabla\phi_m)_{L^2(\omega)}}                                                                                                                  \\
		 & \quad= \biggl|(\chi_\Gamma f_n|\curl \nabla\phi_m)_{L^2(\upsilon)} - (\chi_\Gamma f_n|\curl \nabla\phi_m)_{L^2(\omega)} + \oint_\Gamma \llbracket \nu^\Gamma \times (\chi_\Gamma f_n) \rrbracket \cdot \nabla \phi_m \dd \varsigma\biggr| \\
		 & \quad\leq 0 + \norm{\llbracket \nu^\Gamma \times (\chi_\Gamma f_n) \rrbracket}_{L^2(\Gamma)}\norm{\tr_\tau(\nabla \phi_m)}_{L^2(\Gamma)}                                                                                                  \\
		 & \quad\lesssim \norm{\llbracket \nu^\Gamma \times f_n \rrbracket}_{L^2(\Gamma)}\norm{\phi_m}_{H^1(\Gamma)}.
	\end{align*}
	As $m \to \infty$ we derive
	\[
		\norm{\llbracket\nu^\Gamma \cdot( \curl f_n) \rrbracket}_{H^{-1}(\Gamma)} \lesssim \norm{ \llbracket \nu^\Gamma \times f_n \rrbracket}_{L^2(\Gamma)},
	\]
	and thus the claim letting $n$ tend to infinity.
\end{proof}

\subsection{Div-Curl estimates}

Div-curl estimates allow us to control Sobolev norms through the $L^2$ norms, the divergence and the curl of a given function, if additional information at the boundary is available. Usually this means that either $u \in H_{n0}(\div)$ or $u \in H_{t0}(\curl)$. For example the following version can be found in literature, cf.\ Corollary~3.5 in \cite{amrouche2021}.
\begin{lemma}\label{lem:div-curl-basic}
	Let $k \in \N$, $\Omega \subset \R^3$ with boundary $\partial\Omega$ of class $C^{k,1}$. Assume that $\partial\Omega = \partial\Omega_D \cupdot \partial\Omega_N$ with $\overline{\partial\Omega_D} \cap \overline{\partial\Omega_N} = \emptyset$.
	Let $u \in L^2(\Omega)$ with $\div u \in H^{k - 1}(\Omega)$, $\curl u \in H^{k-1}(\Omega)^3$ such that $u\cdot \nu \in H^{k-1/2}(\partial\Omega_N)$ and $u\times \nu \in H^{k-1/2}(\partial\Omega_D)^3$. Then $u \in H^k(\Omega)$ with
	\begin{align*}
		\norm{u}_{H^k(\Omega)} & \lesssim \norm{u}_{L^2(\Omega)} + \norm{\div u}_{H^{k-1}(\Omega)} + \norm{\curl u}_{H^{k-1}(\Omega)}         \\
		                       & \seq + \norm{u \cdot \nu}_{H^{k-1/2}(\partial\Omega_N)} + \norm{u \times \nu}_{H^{k-1/2}(\partial\Omega_D)}.
	\end{align*}
\end{lemma}

For our analysis, we need a more general variant based upon Proposition~2.24 in \cite{dohnal2024} and Theorem~1.4 in \cite{bresch2026}. We note that the connectedness condition on $\Omega^\complement$ in \cite{dohnal2024} is not needed, see Theorem~1.4 in \cite{bresch2026}.

\begin{theorem}[Div-Curl Lemma]
	\label{thm:div-curl-perfect}
	Let $k \in \N$, $\Omega = \Omega_1 \cupdot \Gamma \cupdot \Omega_2 \subset \R^3$ bounded, where $\Gamma \subset \Omega$ is an interface with $\dist(\Gamma,\partial\Omega)>0$.
	Let $\Gamma$ and $\partial \Omega$ be of class $C^{\hat k+1}$ and $\alpha|_{\Omega_i} \in H^{\hat{k}}(\Omega_i, \Rsym)$, where $\hat{k} \coloneqq \max\{k,2\}$, for
	$i = 1, 2$ be uniformly positive definite, i.e., $\alpha|_{\Omega_i} \geq \eta>0$.
	Let $u \in L^2(\Omega)^3$ satisfy either of	the cases
	\begin{enumerate}
		\item $u \in H_{t0}(\curl, \Omega)$ and $\alpha u + v \in H(\div, \Omega)$, or
		\item $u \in H(\curl, \Omega)$ and $\alpha u + v \in H_{n0}(\div, \Omega)$
	\end{enumerate}
	for some $v \in H_\Gamma^{k}$. Furthermore, let $\curl u \in H_\Gamma^{k-1}(\Omega)$ and $\div(\alpha u)\in H_\Gamma^{k-1}(\Omega)$.
	Then $u$ is contained in $H_\Gamma^{k}$, and it holds
	\begin{equation*}
		\norm{u}_{H^{k}_\Gamma} \le C
		\bigl(\norm{u}_{L^2(\Omega)} + \norm{\curl u}_{H^{k-1}_\Gamma} + \norm{\div (\alpha u)}_{H^{k-1}_\Gamma} + \norm{v}_{H^{k}_\Gamma}\bigr)
	\end{equation*}
	for a constant $C >0$ depending only on $k$, $\smallnorm{\alpha}_{H^{\hat{k}}_\Gamma}$, $\eta$, $\Omega_1$, and $\Omega_2$.
\end{theorem}

\begin{remark}\label{rem:reg-coefficients-div-curl}
	The div-curl lemma is based on regularity results for an elliptic transmission problem, see \cite{bresch2026} and \cite{weber1981}. The cited results are formulated for $C^{k}(\ol{\Omega_\pm})$-coefficients $\alpha$.
	However, it is possible to adapt their proofs to $H^{\hat{k}}_\Gamma$-coefficients, relying on product estimates in Sobolev spaces as in Lemma~2.1 in \cite{schnaubelt2022}.
\end{remark}

\section*{Acknowledgments}
The author thanks Christopher Bresch, David Ploß and Roland Schnaubelt for fruitful discussions and remarks on drafts of this paper.

AI tools, mainly \textit{Gemini 3.1 Pro}, were used to improve the Ti\textit{k}Z code of the figures.
Additionally, the tools were employed in the manipulation and the typesetting in \LaTeX{} of some of the equations present in the text. They were consulted at various points to check results.
An almost finished draft was also given to \mbox{\textit{Claude Fable 5.1}}, which spotted an oversight in one of the proofs, as indicated there.
Other arguments and explanations are human-generated and the author assumes full responsibility for the integrity of the work.

Funded by the Deutsche Forschungsgemeinschaft (DFG, German Research Foundation) -- Project-ID 258734477 -- SFB 1173.

\bigskip

\textbf{Declarations of interest: none.}
\clearpage
\printbibliography 
\end{document}